\documentclass[a4paper, twoside, 11pt]{amsart}
\usepackage{amsthm,amssymb,amsmath,amscd,empheq,mathrsfs}
\usepackage[a4paper, width = 145mm, top = 29mm, bottom = 29mm, bindingoffset= 20mm]{geometry}
\usepackage[shortlabels]{enumitem}
\usepackage{xcolor}
\usepackage[utf8]{inputenc}
\usepackage[foot]{amsaddr}
\usepackage{graphicx,lipsum,caption}
\usepackage{tikz}
\usepackage{verbatim}
\usepackage{tikz-cd}
\usepackage{longtable}
\usepackage{array}
\title{\MakeUppercase{Categorical resolutions of $A_2$ singularities}}
\usepackage{hyperref}
\hypersetup{
  colorlinks   = true,    
  urlcolor     = blue,    
  linkcolor    = blue,    
  citecolor    = teal!80!black      
}

\usepackage[symbol]{footmisc}

\usepackage[style=alphabetic, maxbibnames=99]{biblatex}
\makeatletter
\newtheorem*{rep@theorem}{\rep@title}
\newcommand{\newreptheorem}[2]{%
\newenvironment{rep#1}[1]{%
 \def\rep@title{#2 \ref{##1}}%
 \begin{rep@theorem}}%
 {\end{rep@theorem}}}
\makeatother

\newreptheorem{theorem}{Theorem}

\newtheorem{theorem}{Theorem}[section]
\newtheorem{lemma}[theorem]{Lemma}
\newtheorem{proposition}[theorem]{Proposition}
\newtheorem{corollary}{Corollary}[theorem]

\theoremstyle{definition}
\newtheorem{definition}[theorem]{Definition}
\newtheorem{example}[theorem]{Example}
\theoremstyle{remark} 
\newtheorem{remark}[theorem]{Remark}

\numberwithin{equation}{theorem}

\DeclareMathOperator{\im}{im}

\DeclareMathOperator{\coker}{coker}

\DeclareMathOperator{\Hom}{Hom}

\DeclareMathOperator{\id}{id}

\DeclareMathOperator{\bl}{Bl}
\DeclareMathOperator{\lilperf}{perf}
\DeclareMathOperator{\picard}{Pic}
\DeclareMathOperator{\cone}{cone}
\DeclareMathOperator{\aut}{Aut}

\DeclareMathOperator{\coh}{Coh}
\DeclareMathOperator{\characteris}{char}

\DeclareMathOperator{\codim}{codim}
\DeclareMathOperator{\matri}{M}

\newcommand{\lmod}{{\operatorname{\mathbf{-mod}}}}
\newcommand{\regD}{\widetilde{\mathcal{D}}}
\newcommand{\dbx}{D^b(X)}
\newcommand{\perf}{\mathcal{D}^{\text{perf}}}

\newcommand{\wmax}{W_{\text{max}}}
\newcommand{\evencl}{\mathrm{Cl}_0}
\newcommand{\oddcl}{\mathrm{Cl}_1}
\newcommand{\cl}{\mathrm{Cl}}

\title{Categorical resolutions of Cuspidal singularities}
\author[Céline Fietz]{Céline Fietz \textsuperscript{$\ast$}}
\address{\textsuperscript{$\ast$}Mathematical Institute, Universiteit Leiden, Einsteinweg 55, 2333 CC Leiden, Netherlands.}
\email{c.fietz@math.leidenuniv.nl}
\begin{document}
\maketitle
\begin{abstract}
	Let $X$ be a projective variety of dimension at least two with an isolated $A_2$ singularity. We study its bounded derived category and prove that there exists a crepant categorical resolution $\pi_*\colon \widetilde{\mathcal{D}} \to \dbx$, which is a Verdier localization. More importantly, we give an explicit description of a generating set for its kernel. In the case of an even dimensional variety with a single $A_2$ singularity, we prove that this generating set is given by two $2$-spherical objects. If $X$ is a cubic fourfold with an isolated $A_2$ singularity, we show that this resolution restricts to a crepant categorical resolution $\widetilde{\mathcal{A}}_X$ of the Kuznetsov component $\mathcal{A}_X \subset \dbx$, and that $\widetilde{\mathcal{A}}_X$ is equivalent to the bounded derived category of a K3 surface. 
\end{abstract}

\setcounter{tocdepth}{1} 
\tableofcontents

\reversemarginpar
\section{Introduction}
In recent decades, the study of the bounded derived category of coherent sheaves of a variety has become a powerful and versatile tool in algebraic geometry. While much is known about derived categories of smooth projective varieties, the study of derived categories of singular varieties has become a very active topic only in the last few years. Classically, in order to understand singular varieties, one often studies resolutions of the singularities. This idea has a categorical manifestation, namely the notion of categorical resolutions of a triangulated category. Through this abstraction it is sometimes possible to extend geometric notions that exist only in low dimensions to higher dimensional varieties, such as in the case of simultaneous resolutions of singularities, see \cite{simulataneous}. 

Let $X$ be a projective variety with rational singularities and consider a resolution of singularities $\pi\colon \widetilde{X} \to X$. On the level of bounded derived categories, there exist exact functors  
\begin{equation*}
	\pi_*\colon D^b(\widetilde{X}) \to \dbx \quad \text{and} \quad \pi^*\colon D^{\lilperf}(X) \to D^b(\widetilde{X}),
\end{equation*} 
where $\pi^*$ is left adjoint to $\pi_*$ on $D^{\lilperf}(X)$. Since $X$ has rational singularities, the functor $\pi^*$ is fully faithful. More generally, following \cite[Definition 3.2]{Lefschetz_Decomp}, a categorical resolution of a triangulated category $\mathcal{D}$ is defined as a triple $$(\widetilde{\mathcal{D}}, \pi_*\colon \widetilde{\mathcal{D}} \to \mathcal{D}, \pi^*\colon \perf \to \widetilde{\mathcal{D}}),$$ where the category $\widetilde{\mathcal{D}}$ is a full admissible subcategory of the bounded derived category of a smooth projective variety, there exists an adjunction $\pi^*\vdash \pi_*$ and we require that the functor $\pi^*$ is fully faithful. In the case that $\pi^*$ is additionally right adjoint to $\pi_*$, we call the categorical resolution \textit{crepant}. These resolutions are of special interest, because they are considered to be particularly small\footnote{The definition of crepancy we use in this paper is in some articles referred to as \textit{weak crepancy}. Note that weak crepancy does not imply that the categorical resolution is minimal.}.

Recall that an object $P$ in the bounded derived category $\dbx$ is called $l$\textit{-spherical} if $\Hom(P, P) = k \oplus k[-l]$ and there exists an isomorphism of functors $\Hom(P,-) \cong \Hom(-,P[l])^\vee$. Let $F\colon \mathcal{D} \to \mathcal{D}'$ be an exact functor between triangulated categories. The \textit{kernel of $F$} is defined as the full subcategory $\ker(F)\subset \mathcal{D}$ of objects $A\in \mathcal{D}$ such that $F(A) \cong 0$. We call $F$ a \textit{Verdier localization} if the induced functor $\overline{F}\colon \mathcal{D}/\ker(F) \to \mathcal{D'}$ is an equivalence; here $\mathcal{D}/\ker(F)$ is the \textit{Verdier quotient} of $\mathcal{D}$ by the triangulated subcategory $\ker(F)$.
Our main result is as follows. 
\begin{theorem}\label{MAIN}
   Let $X$ be a projective variety with an isolated $A_2$ singularity (see Definition \ref{defi_A_d}) and $\dim(X) \geq 2$. Then there exists a crepant categorical resolution $\pi_*\colon \widetilde{\mathcal{D}} \to \dbx$ which is a Verdier localization. If $X$ is even dimensional, the kernel $\ker(\pi_*)$
   is generated by two $2$-spherical objects $\mathcal{K}_1, \mathcal{K}_2$. In other words, the functor $\pi_*\colon \widetilde{\mathcal{D}} \to  \dbx$ induces an equivalence of triangulated categories 
	\begin{equation*}
		\overline{\pi}_* \colon \left. \widetilde{\mathcal{D}} \middle/ \langle \mathcal{K}_1,\mathcal{K}_2 \rangle  \right.\overset{\sim}{\longrightarrow} \dbx.
	\end{equation*}
\end{theorem}
In the case of an odd dimensional variety, the kernel is generated by one object $\mathcal{K}$ and we will prove that it is not $l$-spherical for any integer $l$. 

The existence of spherical objects allows us to construct interesting autoequivalences of the category $\widetilde{\mathcal{D}}$. At present it is not known for which kind of singularities on a variety the kernel of a categorical resolution is generated by spherical objects, but if a resolution has this property, it is automatically crepant, see \cite[Lemma 5.8]{hfd}. 

We will show that the existence of the crepant categorical resolution of Theorem \ref{MAIN} follows as a direct application of \cite[§4]{Lefschetz_Decomp}. Moreover, the fact that $\pi_*$ is a Verdier localization will be deduced from a theorem proved by Efimov in \cite[§8]{Efimov_key}. Finally, the main contribution of this paper lies in explicitly determining generators of the kernel $\ker(\pi_*)$ and showing that they are $2$-spherical. 

Note that \cite[Theorem 1.1]{Kernels_nodal}, and simultaneously \cite[Theorem 5.8]{categorical_absorptions}, proved a result similar to Theorem \ref{MAIN} for a variety with an isolated $A_1$ singularity, and our proof follows the same general structure as theirs, essentially because both singularities can be resolved by a single blow-up at the singularity. The crucial difference between $A_1$ and $A_2$ singularities is that their respective exceptional divisors of said blow-up are a smooth or a nodal quadric. Therefore, the spherical objects for $A_1$ singularities are related to spinor bundles on the smooth quadric, and in the $A_2$ case they are related to spinor sheaves on the nodal quadric.

Let us consider the blow-up $\widetilde{X}$ of $X$ at the $A_2$ singularity $x$ which resolves the singularity and provides us with a cartesian diagram 
\begin{equation*}
		\begin{tikzcd}
			Y \arrow[d] \arrow[r, "j", hook] & \widetilde{X} \arrow[d] \\
			\{x\} \arrow[r, hook]            & X.                     
		\end{tikzcd}
	\end{equation*}
Here, $Y$ denotes the exceptional divisor which is a nodal quadric and $j$ denotes the closed embedding $ Y \subset \widetilde{X}$. As we will explain in detail in Subsection \ref{addington_def}, one can define certain reflexive sheaves $\mathcal{S}_1, \mathcal{S}_2$ and $\mathcal{S}$ on $Y$, for $\dim(Y)$ odd and even, respectively, which are locally free away from the singularity, called the \textit{spinor sheaves} of $Y$. Using the notation of Theorem \ref{MAIN}, we will prove that 
\begin{equation*}
    \mathcal{K}_1 = j_*\mathcal{S}_1,  \quad \mathcal{K}_2 = j_*\mathcal{S}_2 \quad \text{and} \quad \mathcal{K} = j_*\mathcal{S}.
\end{equation*}
A key ingredient in proving that $\mathcal{K}_1, \mathcal{K}_2$ are $2$-spherical and that $\mathcal{K}$ is not $l$-spherical for any $l\in \mathbb{Z}$ is the following result:
\begin{theorem}\label{heart_intro}
Let $Y$ be a nodal quadric. If $Y$ is odd dimensional, then there exist isomorphisms of $k$-algebras
        \begin{equation}\label{redo_odd}
            \begin{aligned}
                \mathrm{Ext}^{\bullet}(\mathcal{S}_1,\mathcal{S}_1) \cong \mathrm{Ext}^{\bullet}(\mathcal{S}_2,\mathcal{S}_2)  \cong k[\theta];
                \end{aligned}
        \end{equation}
    where the element $\theta$ has degree $2$. Moreover, the complexes 
        \begin{equation}\label{key}
            \mathrm{Ext}^{\bullet}(\mathcal{S}_1,\mathcal{S}_2) \quad \text{and} \quad \ \mathrm{Ext}^{\bullet}(\mathcal{S}_2,\mathcal{S}_1)
        \end{equation}
    admit a structure of a free $k[\theta]$-module of rank $1$, generated by an element of degree $1$. If $Y$ is even dimensional, there exists an isomorphism of $k$-algebras
        \begin{equation}\label{redo_even}
            \begin{aligned}
                \mathrm{Ext}^{\bullet}(\mathcal{S},\mathcal{S}) \cong k[\theta'],
            \end{aligned}
        \end{equation}
    where the element $\theta'$ has degree $1$.
\end{theorem}
The fact that the self Ext-complexes are isomorphic to polynomial algebras was already proved in \cite[Proposition 6.4]{categorical_absorptions}, and in their terminology, see \cite[Definition 1.7]{categorical_absorptions}, the above result tells us that the spinor sheaves $\mathcal{S}_1, \mathcal{S}_2$ are $\mathbb{P}^{\infty,2}$-objects, and $\mathcal{S}$ is a $\mathbb{P}^{\infty,1}$-object. Our first contribution in Theorem \ref{heart_intro} is to compute the complexes (\ref{key}) which will be needed for showing that $j_*\mathcal{S}_1$ and  $j_*\mathcal{S}_2$ are $2$-spherical. Secondly, we prove Theorem \ref{heart_intro} using the definition of spinor sheaves on nodal quadrics established in \cite{Addington} and thereby give an alternative proof of \cite[Proposition 6.4]{categorical_absorptions}. Although the definitions are equivalent, the one of \cite{Addington} is more algebraic in flavor than the one used in \cite{categorical_absorptions}, and the proof falls out of it naturally. Additionally, the anonymous referee suggested an alternative proof for Theorem \ref{heart_intro} using Kuznetsov--Shinder’s definition. We outline their argument in Subsection 3.4. 

We now give a rough outline of the proof of Theorem \ref{heart_intro} and to this end sketch the construction of the spinor sheaves on a singular quadric, following \cite{Addington}. Let $(V,q)$ be an odd dimensional quadratic space with quadratic form $q$ of corank $1$ and let $Y \subset \mathbb{P}(V)$ denote the associated nodal quadric. We consider the \textit{Clifford algebra} of this quadratic space which is defined as
\begin{equation*}
    \cl(q) := T^{\bullet}(V)/\langle q(v) - v^2 \rangle;
\end{equation*}
here $T^{\bullet}(V)$ denotes the tensor algebra of $V$. This $k$-algebra is in fact $\mathbb{Z}/2\mathbb{Z}$-graded, where $\evencl(q)$ and $\oddcl(q)$ denote the subspaces spanned by elements in $\cl(q)$ that are products of an even and odd number of elements in $V$, respectively. In low dimensions, e.g. for $\dim(V)= 2$ and for $\dim(V) =3$, the Clifford algebra and the even part of the Clifford algebra are given by degenerate quaternions, respectively, see Example \ref{example_intro}. 

We will show that there exists an equivalence of triangulated categories
\begin{equation}
     D^b(\evencl(q)) \overset{\cong}{\longrightarrow} \langle \mathcal{S}_1, \mathcal{S}_2 \rangle \subset D^b(Y),
\end{equation}
where $\langle \mathcal{S}_1, \mathcal{S}_2 \rangle$ denotes the smallest triangulated subcategory in $D^b(Y)$, containing the sheaves $\mathcal{S}_1$ and $\mathcal{S}_2$. In particular, this equivalence shows that we can compute the Ext-complexes of the corresponding $\evencl(q)$-modules. 

The second key step in proving Theorem \ref{heart_intro} is the reduction to a low dimensional quadratic space $(V,q)$ or rather low dimensional Clifford algebras. More precisely, we will show the following version of Knörrer periodicity in Subsection \ref{clifford_alg_Mortia}: 
\begin{proposition}
    Let $(V,q)$ be a quadratic space with $q\neq 0$ and let $U$ denote the hyperbolic plane. Then there exists an isomorphism of $k$-algebras 
    \begin{equation*}
        \evencl(q \perp U) \cong \matri_2(\evencl(q)).
    \end{equation*}
\end{proposition}

Finally, we specialize to a cubic fourfold $X$ with an isolated $A_2$ singularity, whose bounded derived category admits a semiorthogonal decomposition 
\begin{equation*}
    \dbx = \langle \mathcal{A}_X, \mathcal{O}_X,\mathcal{O}_X(1),\mathcal{O}_X(2) \rangle,
\end{equation*}
where $\mathcal{A}_X$ denotes the \textit{Kuznetsov Component}. More precisely, it is defined as the full subcategory 
\begin{equation*}
    \mathcal{A}_X = \{ \mathcal{F} \in \dbx \ | \ \mathrm{Ext}^{\bullet}(\mathcal{O}_X(i), \mathcal{F}) = 0 \ \text{for} \ i = 0,1,2 \ \}.
\end{equation*}
Then the exceptional divisor $Y$ contains a K3 surface $S$, which is deeply related to $\dbx$ as follows.
\begin{theorem}\label{cubic_fourfold_main_result}
	Let $X$ be a cubic fourfold with an isolated $A_2$ singularity and let $\pi\colon \widetilde{X} \to X$ denote the blow-up of $X$ at the singularity $x$. Then there exists a full admissible subcategory $\widetilde{\mathcal{A}}_X \subset D^b(\widetilde{X})$ so that the resolution of Theorem \ref{MAIN} restricts to a crepant categorical resolution 
    \begin{equation*}
        \pi_*\colon\widetilde{\mathcal{A}}_X \longrightarrow \mathcal{A}_X \quad \text{and} \quad \pi^*\colon \mathcal{A}_X^{\lilperf} \to \widetilde{\mathcal{A}}_X,
    \end{equation*}
    of the Kuznetsov component $\mathcal{A}_X$, where $\mathcal{A}_X^{\lilperf} = \mathcal{A}_X \cap D^{\lilperf}(X)$. We have
    \begin{equation*}
        \ker(\widetilde{\mathcal{A}}_X \to \mathcal{A}_X) =  \ker(D^b(\widetilde{X}) \to \dbx) =\langle j_*\mathcal{S}_1, j_*\mathcal{S}_2 \rangle,
    \end{equation*}
    where $j\colon Y \hookrightarrow \widetilde{X}$ denotes the closed immersion of the exceptional divisor $Y$ which is a nodal quadric threefold.
    Moreover, the category $\widetilde{\mathcal{A}}_X$ is equivalent to the derived category $D^b(S)$ of a K3 surface $S$ and under this identification the kernel can be generated by $t^*\mathcal{S}_1$ and $t^*\mathcal{S}_2$, where $t\colon S \hookrightarrow Y$ denotes the embedding of $S$ into the defining nodal quadric $Y$.
\end{theorem}
This generalizes the analogous results for $A_1$ singularities proved in \cite[Theorem 5.2]{dervied_cat_cubic_fourfold} and \cite[§4]{Kernels_nodal}, extending them to $A_2$ singularities with no substantial changes to the arguments. \\
\paragraph{\textbf{Notations and conventions.}}
Throughout the paper we work over a base field $k$ which is algebraically closed and of characteristic $\neq 2$. We call a $k$-scheme $X$ an algebraic variety if it is irreducible, noetherian and separated of finite type over $k$. By $\dbx$ we denote the bounded derived category of coherent sheaves on $X$ which is a $k$-linear triangulated category. The full subcategory of perfect complexes in $X$ is denoted by $\perf(X) \subset \dbx$. Pullback, pushforward, tensor product and $\Hom$ functors are assumed to be derived if not specified otherwise. \\
\paragraph{\textbf{Acknowledgments.}}
This paper is a generalization of my master's thesis which I completed at the University of Bonn in October 2023. I would like to express my deep gratitude to my advisors Yajnaseni Dutta and Evgeny Shinder for the countless discussions and everything they have taught me along the way. I'm very grateful to the anonymous referee for their corrections and comments which greatly improved the structure and readability of the paper, and for sketching an alternative proof for Theorem \ref{heart_intro} which is explained in Subsection \ref{proof_referee}. Moreover, I would also like to thank Nicolas Addington for a discussion about his paper \cite{Addington}. Finally, I thank Greg Andreychev and Omer Bojan for reading several early versions of the paper and helping me to improve the writing. 
\newpage

\section{Preliminaries}\label{prelim}
In this section, we recall the definitions of geometric and categorical resolutions of singularities, along with some fundamental results that we rely on in later sections.
\subsection{Geometric resolutions of singularities}
In this subsection $X,\widetilde{X}$ denote projective $k$-varieties if not specified otherwise. 
\begin{definition}\label{defi_A_d}
	Let $d \geq 1$ and $n = \dim(X)-1$. An isolated singularity $x \in X$ is an \textit{$A_d$ singularity} if there exists an isomorphism
	\begin{equation*}
		\widehat{\mathcal{O}_{X,x}} \cong k[[ x_1,...,x_{n+2}]]/({x_1^2+\dots+x_{n+1}^2+x_{n+2}^{d+1}}).
	\end{equation*}
	In the case $d=1$ we say that $X$ has a \textit{node} at $x$; in the case $d=2$ we say that $X$ has a \textit{cusp} at $x$.
\end{definition}
\begin{definition}
	Let $X$ be normal. Then $X$ has \textit{rational singularities} if we have 
	\begin{equation*}
		R^0\pi_*\mathcal{O}_{\widetilde{X}} \cong \mathcal{O}_X \quad \text{and} \quad R^i\pi_*\mathcal{O}_{\widetilde{X}} = 0 
	\end{equation*}
	for any $i>0$ and every resolution of singularities $\pi\colon \widetilde{X} \to X$, that is,  $\widetilde{X}$ is smooth and $\pi$ is a proper birational morphism. 
\end{definition}
\begin{remark}
    Let X be a projective variety with an isolated $A_d$ singularity and let $\dim(X) \geq 2$. Then $X$ has rational singularities, see \cite{Viehweg_rational}. 
\end{remark}
We consider the following example of a resolution of singularities which will be fundamental for the whole paper.
\begin{lemma}\label{blow_up}
	Let $X$ be a projective variety with an isolated $A_1$ or $A_2$ singularity at a point $x \in X$. Then the blow-up $\pi\colon \widetilde{X} \to X$ of $X$ at the singular point $x$ is a resolution of singularities for $X$. Let $Y\subset \widetilde{X}$ denote the exceptional divisor. Then $Y$ is a smooth quadric if $X$ has an $A_1$ singularity, and $Y$ is a nodal quadric if $X$ has an $A_2$ singularity. Its conormal bundle is given by $\mathcal{N}_{Y/\widetilde{X}}^\vee \cong \mathcal{O}_Y(1)$. Moreover, in both cases, there exists an isomorphism
    \begin{equation}\label{omega_X_tilde}
		\omega_{\widetilde{X}} \cong \pi^*\omega_{X} \otimes \mathcal{O}_{\widetilde{X}}((n-1)Y), 
	\end{equation}
    where $n = \dim(Y) = \dim(X) -1$. 
\end{lemma}
\begin{proof}
    The first part of the statement follows from a standard blow-up computation. Moreover, since $Y\subset \widetilde{X}$ is an exceptional divisor of a blow-up, its conormal is given by $\mathcal{O}_Y(1)$. For (\ref{omega_X_tilde}) we refer to \cite[Proposition 3.5]{Kernels_nodal}, who prove the statement for $A_1$ singularities, and the proof works analogously in the case of $A_2$ singularities.
\end{proof}
\begin{definition}\label{crepant_res}
	Let $\pi\colon \widetilde{X} \to X$ be a resolution of singularities for $X$. Then we say that $\pi$ is \textit{crepant} if there exists an isomorphism $\omega_{\widetilde{X}} \cong \pi^*\omega_X$. 
\end{definition}
\begin{remark}
	In the setting of Lemma \ref{blow_up}, the divisor $(n-1)Y$ on $\widetilde{X}$ is called the discrepancy of $\pi\colon \widetilde{X} \to X$, see \cite[§1.1]{Reid}. This also explains the neologism ``crepant" in the definition above.
\end{remark}
We end the subsection with a result on cubic hypersurfaces with isolated $A_d$ singularities which will be relevant for the proof of Theorem \ref{cubic_fourfold_main_result}. 
\begin{proposition}\label{fundamental}
	Let $X \subseteq \mathbb{P}^{n+2}$ be a cubic hypersurface with an isolated $A_d$ singularity at $x = [1:0:\dots:0]$. Then it is defined by an equation of the form
	\begin{equation*}
		f(x_0,\dots,x_{n+2}) = x_0q(x_1,\dots,x_{n+2}) + g(x_1,\dots, x_{n+2}),
	\end{equation*}
	for suitable quadratic and cubic polynomials $q$ and $g$, respectively. The blow-up $\widetilde{X}$ of $X$ at the singularity $x$ can be identified with the blow-up of $\mathbb{P}^{n+1}$ along the complete intersection $V(q,g)$. Furthermore:
	\begin{itemize}
		\item $d=1$ if and only if $q$ has maximal rank;
		\item $d=2$ implies that $q$ has corank $1$ and $V(g)$ does not pass through the node of $V(q)$.
	\end{itemize}
	In both cases $V(q,g)$ is smooth.
\end{proposition}
\begin{proof}
	Consider the affine neighborhood of $x$, where $x_0 \neq 0$ and the variety $X$ is given by a defining equation of the form
	\begin{equation}\label{defining_equation}
		f(x_1,\dots,x_{n+2}) = c + l(x_1,\dots, x_{n+2})+q(x_1,\dots,x_{n+2}) + g(x_1,\dots, x_{n+2}),
	\end{equation}
	for some homogeneous polynomials $g$, $q$ and $l$ of degrees $3$, $2$ and $1$, respectively, and a scalar $c \in k$. Since the point $0 = (0:\dots:0)$ is a singularity of $f$, we have $l=0$ and $c=0$ which proves the first claim. For the statement about identifying the blow-ups we refer to \cite[§1.5]{geocubic}.  

	Let $d=1$ and consider the blow-up of $X$ at the point $x$. By Lemma \ref{blow_up}, the corresponding exceptional divisor is a smooth quadric, and it can be shown to coincide with the projectivized tangent cone $\mathbb{P}\text{TC}_x(X)$, which is by definition isomorphic to $V(q)$. Conversely, let $q$ be a quadratic form of maximal rank. Since $k$ is algebraically closed and $\characteris(k) \neq 2$ we can assume that $q=x_1^2+\dots +x_{n+2}^2$. Passing to a local analytic neighborhood of the singularity $x \in X$, we can write $f$ as an equation
    \begin{equation*}
        f(x_1,\dots,x_{n+2}) = x_1^2+\dots+x_{n+2}^2, 
    \end{equation*}
    after applying a series of coordinate transformations, see \cite[Theorem 2.46]{Morse_book} and \cite[Morse Lemma (2.2)]{milnorMorse}. 
    
    If $X$ has an $A_2$ singularity at $x$, we can show analogously to the $d=1$ case that $q$ has corank $1$ and therefore assume that $q=x_1^2+\dots +x_{n+1}^2$. To prove the second part of the statement, we start with the following observation: If $g$ contains a term of the form $cx_{n+2}^3$, with $c\in k^{\times}$, then this implies that 
    \begin{equation}\label{observation}
		g(0,0,\dots,0,1) \neq 0
	\end{equation}
    which is equivalent to the fact that the cubic $V(g)$ does not pass through the node of $V(q)$. An elementary computation of the blow-up $\widetilde{X}$ of $X$ at the point $x$ shows that $g$ must contain a term of the form $cx_{n+2}^3$, with $c\in k^{\times}$, because otherwise $\widetilde{X}$ would be singular, which is a contradiction.

    Lastly, the smoothness of $V(q,g)$ can be proven in an elementary way, like in \cite[Proposition 2.8]{A2_v2}, or by the fact that we know that $\widetilde{X}$ is smooth for $A_1$ and $A_2$ singularities; therefore we can apply 
    \cite[Lemma 2.6]{blow_up_smooth_center}. 
\end{proof}
\begin{corollary}\label{S_is_K3}
	Let $X \subset \mathbb{P}^5$ be a cubic fourfold which is smooth away from an isolated $A_2$ singularity at a point $x \in X$. Then, with the notation of Proposition \ref{fundamental}, the intersection $S = V(q,g)$ is a smooth K3 surface. 
\end{corollary}
\begin{proof}
	In Proposition \ref{fundamental} we proved the smoothness of $S = V(q,g)$, so it remains to show that $S$ is indeed a K3 surface. 
	The adjunction formula gives us an isomorphism $\omega_S \cong \mathcal{O}_S(2+3-5) = \mathcal{O}_S$, and since $S$ is a complete intersection, a standard calculation shows that $H^1(S, \mathcal{O}_S) = 0$, see for example \cite[Ex.III.5.5c]{hartshorne}. 
\end{proof}
\subsection{Semiorthogonal decompositions, Serre functors and mutations}
Let $\mathcal{D}$ be a triangulated category and $\mathcal{A} \subset \mathcal{D}$ a full triangulated subcategory. The \textit{left orthogonal} to  $\mathcal{A}$ in $\mathcal{D}$ is defined as the full triangulated subcategory  
\begin{equation*}
	^\perp\mathcal{A} = \{B \in \mathcal{D} \ | \ \Hom_{\mathcal{D}}(B,A) = 0 \ \forall A \in \mathcal{A}\}.
\end{equation*}
Analogously, the \textit{right orthogonal} to $\mathcal{A}$ in $\mathcal{D}$ is defined as the full triangulated subcategory  
\begin{equation*}
	\mathcal{A}^\perp =\{B \in \mathcal{D} \ | \ \Hom_{\mathcal{D}}(A,B) = 0 \ \forall A \in \mathcal{A}\}.
\end{equation*}
\begin{definition}[{\cite{Bondal_ass_alg_coh}, \cite{BK_rep_funcots_serre_func}}]
	Let $\mathcal{A} \subset \mathcal{D}$ denote a full triangulated subcategory and denote the inclusion functor by $i_*\colon \mathcal{A} \to \mathcal{D}$. Then $\mathcal{A}$ is called \textit{left admissible} (resp.\ \textit{right admissible}) if $i_*$ admits a left adjoint $i^*\colon \mathcal{D} \to \mathcal{A}$ (resp.\ a right adjoint $i^!\colon \mathcal{D} \to \mathcal{A}$). If $\mathcal{A}$ is left and right admissible, then we call it an \textit{admissible} subcategory of $\mathcal{D}$.
\end{definition}
\begin{definition}(\cite{BK_rep_funcots_serre_func}, \cite{Bondal_Orlov_SOD})\label{SOD_def}
	A \textit{semiorthogonal decomposition} of $\mathcal{D}$ consists of full triangulated subcategories $\mathcal{A}_1, \dots, \mathcal{A}_n$, such that 
	\begin{enumerate}
		\item the sequence $\mathcal{A}_1, \dots, \mathcal{A}_n$ is \textit{semiorthogonal}, i.e.\,
		\begin{equation*}
			\Hom_{\mathcal{D}}(A_j, A_i) = 0
		\end{equation*} 
		for all $j >i$ and $A_j\in \mathcal{A}_j,A_i \in \mathcal{A}_i$.
		\item The category $\mathcal{D}$ is the smallest triangulated subcategory of $\mathcal{D}$ containing the subcategories $\mathcal{A}_1, \dots, \mathcal{A}_n$.
	\end{enumerate}
We denote a semiorthogonal decomposition by $\mathcal{D} = \langle \mathcal{A}_1, \mathcal{A}_2, \dots, \mathcal{A}_n \rangle$. We call it \textit{admissible} if all the subcategories $\mathcal{A}_1, \dots, \mathcal{A}_n$ are admissible. 
\begin{remark}\label{admissible_if_sm}
	Let $X$ be a smooth projective variety. Then any semiorthogonal decomposition 
	$\dbx= \langle \mathcal{A}_1, \dots, \mathcal{A}_n \rangle$ is admissible, see \cite{BK_rep_funcots_serre_func}.
\end{remark}
\end{definition}
\begin{lemma}[{\cite[Lemma 3.1]{Bondal_ass_alg_coh}}]\label{SOD}
	Let $\mathcal{A}_1, \mathcal{A}_2, \dots, \mathcal{A}_n$ be a semiorthogonal sequence in $\mathcal{D}$, such that $\mathcal{A}_1, \dots, \mathcal{A}_k$ are \textit{left admissible} and $\mathcal{A}_{k+1}, \dots, \mathcal{A}_n$ are \textit{right admissible}, then 
	\begin{equation*}
		\mathcal{D} = \langle \mathcal{A}_1, \dots, \mathcal{A}_k, ^\bot\langle \mathcal{A}_1, \dots, \mathcal{A}_k \rangle \cap \langle \mathcal{A}_{k+1}, \dots, \mathcal{A}_n \rangle^\bot, \mathcal{A}_{k+1}, \dots, \mathcal{A}_n \rangle
	\end{equation*}
is a semiorthogonal decomposition.
\end{lemma}
It is a well-known fact that one can associate to 
any hypersurface $X \subset \mathbb{P}^{n+1}$ an admissible subcategory in $D^b(X)$ which, intuitively, captures all the ``interesting'' geometric information about $X$. More precisely:
\begin{proposition}\label{Kuzcomp_for_gorenstein}
	Let $X \subset \mathbb{P}^{n+1}$ be a hypersurface of degree $d$ and assume $d \leq n+1$. Then $(\mathcal{O}_X, \dots, \mathcal{O}_X(n+1-d))$ is an exceptional collection in $D^b(X)$. Let 
	\begin{equation*}
		\mathcal{A}_X := \langle \mathcal{O}_X, \dots, \mathcal{O}_X(n+1-d) \rangle^\perp
	\end{equation*}
	denote the right orthogonal of this collection. Then we have
	\begin{equation*}
		\dbx = \langle \mathcal{A}_X, \mathcal{O}_X, \dots, \mathcal{O}_X(n+1-d) \rangle;
	\end{equation*}
    and we will call the subcategory $\mathcal{A}_X$ the \textit{Kuznetsov Component} of $D^b(X)$.
\end{proposition}
\begin{definition}[{\cite[Definition 3.1]{BK_rep_funcots_serre_func}}]
	Let $\mathcal{D}$ be a $k$-linear Hom-finite triangulated category, that is, $\Hom(A,B)$ is finite dimensional for any objects $A,B \in \mathcal{D}$. Then an exact equivalence $\mathbb{S}_{\mathcal{D}} \colon\mathcal{D} \to \mathcal{D}$ is called \textit{Serre functor} if there exist natural isomorphisms 
	\begin{equation*}
		\Hom(F,G) \cong \Hom(G, \mathbb{S}_{\mathcal{D}}(F))^\vee
	\end{equation*}
	of $k$-vector spaces for any $F,G \in \mathcal{D}$. 
\end{definition}
\begin{example}[{\cite[Example 3.2]{BK_rep_funcots_serre_func}}]\label{Serre_smooth}
	Let $X$ be a smooth projective variety. Then $\dbx$ admits a Serre functor which is given by $(-) \otimes \omega_{X}[\dim(X)]$.
\end{example}
\begin{example}[{\cite[§2.3]{Kalck_Shinder_Pavic}, see also \cite[§6]{Ballard}}]\label{GVD_Gorenstein}
	Let $X$ be a Gorenstein projective variety. Then the restriction of $\mathbb{S}_X$ to $D^{\lilperf}(X)$ defines a Serre functor on $D^{\lilperf}(X)$. More generally, there exist natural isomorphisms 
	\begin{equation*}
		\Hom({F}, {G}) \cong \Hom({G},\mathbb{S}_X(F))^\vee,
	\end{equation*}
	for all ${F} \in D^{\lilperf}(X)$ and ${G} \in \dbx$, by Grothendieck-Verdier duality. 
\end{example}
For any triangulated category the group of autoequivalences $\aut(\mathcal{D})$ naturally acts on the set of semiorthogonal decompositions of $\mathcal{D}$. In the special case of Serre functors we recall the following result.
\begin{lemma}[{\cite[Lemma 2.15]{Kalck_Shinder_Pavic}}]\label{SOD_and_Serre}
	Let $X$ be a Gorenstein projective variety and let $\mathcal{D} = \langle \mathcal{A}, \mathcal{B} \rangle$ be a semiorthogonal decomposition, where either $\mathcal{A}$ or $\mathcal{B}$ is contained in $D^{\lilperf}(X)$. Then the subcategories $\mathcal{A}, \mathcal{B} \subset \dbx$ are admissible and there exist semiorthogonal decompositions 
	\begin{equation*}
		\dbx = \langle \mathcal{B}\otimes \omega_{X}, \mathcal{A} \rangle = \langle \mathcal{B}, \mathcal{A} \otimes \omega_X^\vee \rangle.
	\end{equation*}
\end{lemma}
Apart from the automorphisms induced from elements of $\aut(\mathcal{D})$, there exists an interesting class of functors, the \textit{mutation functors} which, informally speaking, permute the components of a semiorthogonal decomposition. 
\begin{proposition}[{\cite{Bondal_ass_alg_coh}}]
	Let $\mathcal{A} \subset \mathcal{D}$ be an admissible subcategory. By Lemma \ref{SOD} we have semiorthogonal decompositions $\mathcal{D} =\langle \mathcal{A}^\bot, \mathcal{A}\rangle$ and $\mathcal{D} =\langle \mathcal{A}, ^\bot\mathcal{A}\rangle$. Then there exist functors $\mathbb{L}_{\mathcal{A}}, \mathbb{R}_{\mathcal{A}}\colon \mathcal{D} \to \mathcal{D}$, vanishing on $\mathcal{A}$ that restrict to mutually inverse equivalences $\mathbb{L}_{\mathcal{A}}\colon^\bot\mathcal{A} \to \mathcal{A}^\bot$ and $\mathbb{R}_{\mathcal{A}}\colon\mathcal{A}^\bot \to \ ^\bot\mathcal{A}$.
\end{proposition}
	We call the functors $\mathbb{L}_{\mathcal{A}}$ and $\mathbb{R}_\mathcal{A}$ the \textit{left-} and \textit{right mutation functors} corresponding to $\mathcal{A}$. 
\begin{lemma}[{\cite{Bondal_ass_alg_coh}}]
	Let $\mathcal{D} = \langle \mathcal{A}_1, \mathcal{A}_2, \dots, \mathcal{A}_n\rangle$ be an admissible semiorthogonal decomposition. Then for each $1\leq k \leq n-1$, there exists a semiorthogonal decomposition
	\begin{equation*}
		\mathcal{D} = \langle \mathcal{A}_1, \dots,  \mathcal{A}_{k-1}, \mathbb{L}_{\mathcal{A}_k}(\mathcal{A}_{k+1}), \mathcal{A}_k, \mathcal{A}_{k+2}, \dots,\mathcal{A}_n\rangle.
	\end{equation*}
	Furthermore, for each $2 \leq k \leq n$, there exists a semiorthogonal decomposition 
	\begin{equation*}
		\mathcal{D} = \langle \mathcal{A}_1, \dots, \mathcal{A}_{k-2}, \mathcal{A}_k, \mathbb{R}_{\mathcal{A}_k}(\mathcal{A}_{k-1}), \mathcal{A}_{k+1}, \dots , \mathcal{A}_n\rangle.
	\end{equation*}
\end{lemma}
\subsection{Categorical resolutions of singularities}
We recall the definition of, and a method for constructing, (crepant) categorical resolutions of the bounded derived category $\dbx$ of a variety $X$ with rational singularities, following \cite{Lefschetz_Decomp}. 
\begin{definition}\label{smooth_triang}
	Let $\mathcal{D}$ be a triangulated category. We say that $\mathcal{D}$ is \textit{smooth} if there exists a smooth projective variety $\widetilde{X}$ such that $\mathcal{D}$ is equivalent to a full admissible subcategory of the bounded derived category $D^b(\widetilde{X})$.
\end{definition}
\begin{remark}
	Nowadays, Definition \ref{smooth_triang} is considered outdated, but we will still use it as it is sufficient for our purposes. We refer to \cite[Definition 3.23]{Orlov_smooth} for a more general definition of smoothness for any enhanced triangulated category. 
\end{remark}
\begin{definition}[{\cite[Definition 1.6]{Orlov_perf}}]\label{defi_hom_finite}
	Let $\mathcal{D}$ be a $k$-linear triangulated category. An object $F \in \mathcal{D}$ is said to be \textit{homologically finite} if for any $G \in \mathcal{D}$ there is only a finite number of values $i \in \mathbb{Z}$, such that  $\Hom_{\mathcal{D}}(F,G[i]) \neq 0$. We denote the full triangulated subcategory of homologically finite objects by $\mathcal{D}^{\lilperf}$. 
\end{definition}
\begin{remark}[{\cite[Lemma 1.11]{Orlov_perf}}]
	Let $X$ be a quasi-projective variety. For $\mathcal{D} = \dbx$ the subcategory of perfect complexes coincides with the subcategory of homologically finite objects which justifies the notation $\mathcal{D}^{\lilperf}$. 
\end{remark}
\begin{definition}[{\cite[Definition 3.2]{Lefschetz_Decomp}}] \label{catresdefi}
	A \textit{categorical resolution} of a triangulated category $\mathcal{D}$ consists of a smooth (as in Definition \ref{smooth_triang}) triangulated category $\regD$ and a pair of functors 
	\begin{equation*}
		\pi_*\colon \regD \to \mathcal{D} \quad \text{and} \quad \pi^*\colon \perf \to \regD
	\end{equation*}
	satisfying the following properties:
	\begin{enumerate}
		\item The functor $\pi^*$ is left adjoint to $\pi_*$. That is, there exist natural isomorphisms
		\begin{equation*}
			\Hom_{\regD}(\pi^*F, G) \cong \Hom_{\mathcal{D}}(F,\pi_*G) \quad \quad \text{for any} \ F\in \perf, \ G \in \regD. 
		\end{equation*}
		\item The natural transformation $\id_{\perf} \to \pi_*\pi^*$ is an isomorphism.
	\end{enumerate}
\end{definition}
\begin{remark}
	Let $X$ be a variety with rational singularities and $\pi\colon \widetilde{X} \to X$ a resolution of singularities. Then the derived category $D^b(\widetilde{X})$ together with the pushforward and pullback functors $\pi_*\colon D^b(\widetilde{X}) \to \dbx$ and $\pi^*\colon D^{\lilperf}(X) \to D^b(\widetilde{X})$ is a categorical resolution of $\dbx$. 
	By imposing the second condition in the above definition, we restrict ourselves to the case where $X$ has at most rational singularities.
\end{remark}
\begin{definition}[{\cite[Definition 3.4]{Lefschetz_Decomp}}] \label{crepcatresdefi}
	A categorical resolution $(\widetilde{\mathcal{D}}, \pi_*,\pi^*)$ is called  \textit{crepant} if the functor $\pi^*$ is also right adjoint to $\pi_*$ when restricted to $\perf$. In other words, there exist natural isomorphisms
	\begin{equation*}
		\Hom_{\regD}(G, \pi^*F) \cong \Hom_{\mathcal{D}}(\pi_*G, F) \quad \quad \text{for any} \ F\in \perf, \ G \in \regD.
	\end{equation*}
\end{definition}
\begin{remark}\label{crepancy_GVD}
	Let $X$ be a Gorenstein projective variety. A  crepant (geometric) resolution $\pi\colon \widetilde{X} \to X$ induces a crepant categorical resolution $(D^b(\widetilde{X}), \pi_*, \pi^*)$ of $\dbx$ as follows. The relative canonical complex $\omega_{\pi}$ is in fact a line bundle, since $X$ is Gorenstein and by Grothendieck-Verdier duality the right adjoint $\pi^!$ of $\pi_*\colon D^b(\widetilde{X}) \to \dbx$ can be given explicitly by $\pi^!(F) = \pi^*F \otimes \omega_\pi$ for all $F \in D^{\lilperf}(X)$. If $\pi$ is crepant then $\omega_\pi \cong \mathcal{O}_{\widetilde{X}}$, hence $\pi^! = \pi^*$.
\end{remark}
In order to describe a construction of a (crepant) categorical resolution, we recall the following definition.
\begin{definition}[{\cite[Definition 2.16]{Lefschetz_Decomp}}]\label{Dual_Lefschetz_decomp}
	Let $X$ be a variety and let $\mathcal{O}(1)$ denote a line bundle on $X$. A \textit{Lefschetz decomposition} of $\dbx$  is a semiorthogonal decomposition of the form 
	\begin{equation*}
		\dbx = \langle \mathcal{B}_0, \mathcal{B}_1(1), \dots,\mathcal{B}_{m-1}(m-1) \rangle, 
	\end{equation*}
	where $\mathcal{B}_0, \mathcal{B}_1, \dots, \mathcal{B}_{m-1}$ are subcategories of $\dbx$, satisfying
	\begin{equation*}
		0 \subset \mathcal{B}_{m-1} \subset \dots \subset \mathcal{B}_1 \subset \mathcal{B}_0 \subset \dbx.
	\end{equation*}
	Similarly, a \textit{dual Lefschetz decomposition} of $\dbx$  is a semiorthogonal decomposition of the form 
	\begin{align*}
	    \dbx = \langle \mathcal{B}_{m-1}(1-m),\dots, \mathcal{B}_1(-1), \mathcal{B}_0\rangle, \  \text{where} \ \ 0 \subset \mathcal{B}_{m-1} \subset \dots \subset \mathcal{B}_1 \subset \mathcal{B}_0 \subset \dbx. 
	\end{align*}
\end{definition}
Following \cite[§4]{Lefschetz_Decomp}, let $X$ denote a projective variety with at most rational singularities and $\pi\colon \widetilde{X} \to X$ a geometric resolution, for which we assume that the exceptional locus $E$ is an irreducible divisor. Let $Z$ be the image of $E$ under $\pi$. Then we have a cartesian diagram 
\begin{equation}\label{gerneral_square}
	\begin{tikzcd}
		E \arrow[d, "p"'] \arrow[r, "j"] & \widetilde{X} \arrow[d, "\pi"] \\
		Z \arrow[r, "i"]                     & X,                             
	\end{tikzcd}
\end{equation}
where the morphisms $i$ and $j$ denote the respective inclusions of the subvarieties $Z \subset X$ and $E \subset \widetilde{X}$.
\begin{theorem}[{\cite[Lemma 4.1, Theorem 4.4]{Lefschetz_Decomp}}]\label{LefschetzProp}\label{catres_withLefschetz}
    In the setting of the preceding paragraph assume that there exists a dual Lefschetz decomposition
	\begin{equation}\label{dual_Lefschetz}
		D^b(E) = \langle \mathcal{B}_{m-1}(1-m), \mathcal{B}_{m-2}(2-m),\dots, \mathcal{B}_1(-1), \mathcal{B}_0\rangle 
	\end{equation} 
	with respect to the conormal bundle $\mathcal{O}_E(1) := \mathcal{N}_{E/\widetilde{X}}^{\vee}$. We define a full triangulated subcategory of $D^b(\widetilde{X})$ by
	\begin{equation}\label{D_tilde_defi}
		\widetilde{\mathcal{D}} = \{F \in D^b(\widetilde{X}) \ | \ j^*F \in \mathcal{B}_0\}.
	\end{equation}
	Then the functor $j_*\colon D^b(E) \to D^b(\widetilde{X})$ is fully faithful when restricted to the subcategories $\langle \mathcal{B}_k(-k) \rangle$ for all $1 \leq k \leq m-1$ and there exists a semiorthogonal decomposition 
	\begin{equation*}
		D^b(\widetilde{X}) = \langle j_*\mathcal{B}_{m-1}(1-m), j_*\mathcal{B}_{m-2}(2-m), \dots, j_*\mathcal{B}_{1}(-1),\widetilde{\mathcal{D}}\rangle.
	\end{equation*}
    If the image of the pullback functor $\pi^*\colon D^{\lilperf}(X) \to D^b(\widetilde{X})$ is contained in $\widetilde{\mathcal{D}}$, the triple $(\widetilde{\mathcal{D}}, \pi_*, \pi^*)$ is a categorical resolution of $\dbx$. 
\end{theorem}
\begin{lemma}[{\cite[Lemma 4.3]{Lefschetz_Decomp}}]\label{small_to_big}
    Assume the setting of the preceding theorem. Then the image of the pullback functor $\pi^*\colon D^{\lilperf}(X) \to D^b(\widetilde{X})$ is contained in $\widetilde{\mathcal{D}}$ if the image of the pullback functor $p^*\colon D^{\lilperf}(Z) \to D^b(E)$ is contained in $\mathcal{B}_0$.
\end{lemma}
\begin{proof}
    Let $F\in D^{\lilperf}(X)$ and consider the diagram (\ref{gerneral_square}). Then the object
    \begin{equation*}
        j^*\pi^*F = p^*i^*F 
    \end{equation*}
    is contained in $\mathcal{B}_0$ by assumption. By definition (\ref{D_tilde_defi}),
    we conclude that $\pi^*F \in \widetilde{\mathcal{D}}$.
\end{proof}
The following result gives sufficient conditions under which the categorical resolution of the previous theorem is crepant. 
\begin{proposition}[{\cite[Proposition 4.5]{Lefschetz_Decomp}}]\label{crepancy_check}
	Let $X$ be Gorenstein and assume that we have an inclusion $p^*(D^{\lilperf}(Z)) \subset \mathcal{B}_{m-1}$. Furthermore, we assume that there exists an isomorphism $\omega_{\widetilde{X}} \cong \pi^*\omega_X \otimes \mathcal{O}_{\widetilde{X}}((m-1)E)$. Then the categorical resolution $(\widetilde{\mathcal{D}}, \pi_*,\pi^*)$ is crepant. 
\end{proposition}

\section{Spinor sheaves on singular quadrics and Morita reduction}\label{spinors_on_singular_quadrics}
The aim of this section is to prove Theorem \ref{heart_intro}, using the definition of spinor sheaves on singular quadrics introduced in \cite{Addington}. Throughout the section $(V,q)$ denotes an arbitrary nontrivial quadratic space, meaning $V \neq \{0\}$. We set $N := \dim(V)$, and therefore the associated quadric hypersurface $Y :=V(q) \subset \mathbb{P}(V)$ has dimension $n:=N-2$.
\subsection{Spinor sheaves on singular quadrics}\label{addington_def} The Clifford algebra $\cl(q)$ associated to $(V,q)$ is defined as the quotient 
\begin{equation}
    \cl(q) := T^{\bullet}(V)/\langle q(v) - v^2 \rangle,
\end{equation}
where $T^{\bullet}(V)$ denotes the tensor algebra. This algebra carries a natural $\mathbb{Z}/2\mathbb{Z}$-grading induced by the involution $V \to V, \ v \mapsto -v$. We denote the eigenspaces with respect to this action by $\evencl(q)$ and $\oddcl(q)$, and obtain a grading $\cl(q) = \evencl(q) \oplus \oddcl(q)$. 

Let $W \subset V$ be an isotropic subspace and choose a basis $\{w_1,\dots, w_m\}$ of $W$. Then one can define a left ideal $I^W = \cl(q)w_1 \cdots w_m \subset \cl(q)$ which admits an induced $\mathbb{Z}/2\mathbb{Z}$-grading $I^W = I^W_{0} \oplus I^W_{1}$. Note that $I^W$ is independent of the choice of basis of $W$, see \cite[§2]{Addington}. The grading gives rise to $k$-linear maps 
\begin{equation}\label{multiplication}
    I^W_{0} \overset{v\cdot}{\longrightarrow} I^W_{1}, \quad \quad I^W_{1} \overset{v\cdot}{\longrightarrow} I^W_{0}
\end{equation}
defined by left multiplication with a vector $v\in V$. Let $\{v_1,\dots,v_{N}\}$ denote a $k$-basis for the vector space $V$ 
and let $x_1,\dots,x_{N}$ be the coordinates of $\mathbb{P}(V)$; then the embedding 
\begin{equation*}
    \mathcal{O}_{\mathbb{P}(V)} \hookrightarrow \mathcal{O}_{\mathbb{P}(V)}(1) \otimes V, \quad 1 \mapsto \sum_i x_i \otimes v_i
\end{equation*} 
induces a morphism
\begin{equation*}
\begin{aligned}
    \mathcal{O}\otimes I_0^W \longrightarrow \mathcal{O}(1) \otimes V \otimes I_0^W, \quad 1\otimes w \longmapsto (\sum_ix_i\otimes v_i) \otimes w. 
\end{aligned}
\end{equation*}
By composition with the multiplication map, we obtain a morphism 
\begin{equation*}
\begin{aligned}
    \varphi\colon \mathcal{O}\otimes I_0^W &\longrightarrow \mathcal{O}(1) \otimes V \otimes I_0^W  \longrightarrow \mathcal{O}(1) \otimes I_1^W,  \\
    1\otimes w &\longmapsto (\sum_ix_i\otimes v_i) \otimes w \longmapsto \sum_i x_i \otimes v_iw
\end{aligned}
\end{equation*}
of vector bundles on $\mathbb{P}(V)$. Analogously, by swapping the roles of $I_0^W$ and $I_1^W$, we define a morphism of vector bundles
\begin{equation*}
\begin{aligned}
    \psi\colon \mathcal{O}\otimes I_1^W \longrightarrow \mathcal{O}(1) \otimes I_0^W,  \\
    1\otimes w' \longmapsto \sum_i x_i \otimes v_iw'.
\end{aligned}
\end{equation*}
\begin{lemma}[{\cite[§2]{Addington}}]\label{phi_and_psi_calc}
    Let $W\subset V$ be an isotropic subspace. Then the compositions 
    \begin{equation*}
    \begin{aligned}
    \mathcal{O} \otimes I_{0}^W \overset{\varphi}{\longrightarrow} \mathcal{O}(1) \otimes I_{1}^W \overset{\psi}{\longrightarrow} \mathcal{O}(2) \otimes I_{0}^W; \\
    \mathcal{O} \otimes I_{1}^W \overset{\psi}{\longrightarrow} \mathcal{O}(1) \otimes I_{0}^W \overset{\varphi}{\longrightarrow} \mathcal{O}(2) \otimes I_{1}^W,
    \end{aligned} 
    \end{equation*}
    coincide with the maps
    \begin{equation*}
    \begin{aligned}
         q \otimes \id_{I_0} \colon \mathcal{O} \otimes I_{0}^W \longrightarrow \mathcal{O}(2) \otimes I_{0}^W,\quad 1 \otimes w \mapsto q \otimes w; \\
         q \otimes \id_{I_1} \colon \mathcal{O} \otimes I_{1}^W \longrightarrow \mathcal{O}(2) \otimes I_{1}^W,\quad 1 \otimes w' \mapsto q \otimes w',
    \end{aligned}
    \end{equation*}
    respectively. In particular, the morphisms $\varphi$ and $\psi$ are injective, and furthermore isomorphisms on the locus where $q\neq 0$.
\end{lemma}
The fact that the cokernels of $\varphi$ and $\psi$ are supported on the quadric $Y$ motivates the following definition.
\begin{definition}[{\cite[§2]{Addington}}]\label{definition_spinor}
    Let $W \subset V$ be an isotropic subspace, then we call the sheaves 
    \begin{equation*}
        \mathcal{S}^W :=\coker(\varphi(-1)), \ \ \text{and} \ \ \mathcal{T}^W:= \coker(\psi(-1))
    \end{equation*}
    on $Y$ the \textit{spinor sheaves associated to $W$}. Since $\varphi$ and $\psi$ are injective, we have the following short exact sequences of sheaves 
    \begin{equation}\label{defining_sequences_spinor}
    \begin{aligned}
        0 \longrightarrow \mathcal{O}_{\mathbb{P}(V)}(-2)^M\overset{\varphi}{\longrightarrow} \mathcal{O}_{\mathbb{P}(V)}^M(-1)\longrightarrow \mathcal{S}^W \longrightarrow 0; \\
        0 \longrightarrow \mathcal{O}_{\mathbb{P}(V)}(-2)^M\overset{\psi}{\longrightarrow} \mathcal{O}_{\mathbb{P}(V)}^M(-1) \longrightarrow \mathcal{T}^W \longrightarrow 0,
    \end{aligned}
    \end{equation} 
    where $M = \dim(I_0) = \dim(I_1) = 2^{\codim(W)-1}$.
\end{definition}
\begin{remark}
    In \cite[§2]{Addington} the spinor sheaves are defined by the same short exact sequences (\ref{defining_sequences_spinor}), but are additionally twisted by $\mathcal{O}_{\mathbb{P}(V)}(1)$. This convention implies that the spinor bundle on a smooth conic $Q$ is given by $\mathcal{O}_Q(1)$ and coincides with \cite{Kapranov_spinors}. In contrary, using the above convention, the spinor bundle on $Q$ will coincide with $\mathcal{O}_Q(-1)$, as in \cite{Ottaviani_spinors}.
\end{remark}
With the short exact sequences (\ref{defining_sequences_spinor}) above, we can easily deduce the following result about the cohomology of spinor sheaves. 
\begin{lemma}\label{cohomology}
    Let $(V,q)$ be a quadratic space with $\dim(V) \geq 3$ and let $Y \subset \mathbb{P}(V)$ denote the associated quadric. Then for any isotropic subspace $W\subset V$ the cohomology of the spinor sheaves $\mathcal{S}^W, \mathcal{T}^W$ vanishes. Moreover, we have 
    \begin{equation*}
        H^i(Y,\mathcal{S}^W(l)) = 0 \quad \text{and} \quad  H^i(Y,\mathcal{T}^W(l)) = 0,
    \end{equation*}
    for $0<i<n$ and any $l \in \mathbb{Z}$, that is, their intermediate cohomology vanishes.
\end{lemma}
Furthermore, there exist the following short exact sequences that will be used frequently in the rest of the paper. 
\begin{lemma}\label{exactsequences_spinors_o}
    Let $(V,q)$ be a quadratic space, let $Y \subset \mathbb{P}(V)$ denote the associated quadric and let $W \subset V$ be an isotropic subspace. Then there exist short exact sequences 
    \begin{equation}\label{ses_with_twist}
        \begin{aligned}
            0 \longrightarrow \mathcal{T}^W\overset{\varphi}{\longrightarrow} \mathcal{O}_Y^M \longrightarrow \mathcal{S}^W(1)\longrightarrow 0; \\
            0 \longrightarrow \mathcal{S}^W \overset{\psi}{\longrightarrow} \mathcal{O}_Y^M\longrightarrow \mathcal{T}^W(1) \longrightarrow 0,
        \end{aligned}
    \end{equation}
    where $M = \dim_k(I_0^W) = \dim_k(I_1^W) = 2^{\codim(W)-1}$.
\end{lemma}
\begin{proof}
    The short exact sequences can be constructed by breaking up the locally free resolutions in \cite[§4]{Addington}.
\end{proof}
Given a quadratic space $(V,q)$, we will denote the kernel of the quadratic form $q\colon V \to k$ by $K$. Then the singular locus of the quadric $Y =V(q) \subset \mathbb{P}(V)$ coincides with $\mathbb{P}(K) \subset \mathbb{P}(V)$.
\begin{proposition}[{\cite[Proposition 2.1]{Addington}}]\label{locallyfree_or_not}
    The restrictions of the spinor sheaves $\mathcal{S}^W,\mathcal{T}^W$ to $\mathbb{P}(K) \cap \mathbb{P}(W)$ are trivial of rank $2^{\codim(W)-1}$. In the case that $\codim(W) > 1$, the sheaves $\mathcal{S}^W$ and $\mathcal{T}^W$ are locally free of rank $2^{\codim(W)-2}$ on $Y\setminus (\mathbb{P}(K) \cap \mathbb{P}(W)) $. In particular, $\mathcal{S}^W$ and $\mathcal{T}^W$ are locally free if and only if $K\cap W =\{0\}$.
\end{proposition}
\begin{remark}[{\cite[§3]{Addington}}]
    The sheaves $\mathcal{S}^W$, $\mathcal{T}^W$ are unchanged while varying $W$ continuously with $W \cap K$ fixed. Moreover, let $\pi\colon V \to V/K$ be the projection, then the following holds:
    \begin{enumerate}[label=\roman*)]
        \item If $\dim(\pi(W)) < \frac{1}{2} \dim(V/K)$, then $\mathcal{S}^W \cong \mathcal{T}^W$;
        \item If $\dim(\pi(W)) = \frac{1}{2}\dim(V/K)$, then $\mathcal{S}^W \not\cong \mathcal{T}^W$ and switching $\pi(W)$ to the other connected family (while keeping $W\cap K$ fixed) interchanges $\mathcal{S}^W$ and $\mathcal{T}^W$. 
    \end{enumerate}
\end{remark}
\begin{definition}\label{defi_spinors_sec3}
    Let $(V,q)$ be a quadratic space with $q$ of corank $1$. We denote the spinor sheaves corresponding to the maximal isotropic subspace $\wmax \subset V$, for $V$ odd and even dimensional, by $\mathcal{S}_1,\mathcal{S}_2$ and $\mathcal{S}$, respectively. We will refer to them as \textit{``the" spinor sheaves} on the nodal quadric $Y = V(q)$.
\end{definition}
\begin{remark}
    In general, the dimension of the maximal isotropic subspace $\wmax \subset V$ is given by $$d_{\mathrm{max}} = \left\lfloor\frac{\dim(V/K)}{2} \right\rfloor +\dim(K)= \left\lfloor\frac{\dim(V)+\dim(K)}{2}\right\rfloor.$$
\end{remark}
In the remainder of this subsection, we classify all spinor sheaves on a nodal quadric $Y$, up to isomorphism. The key ingredient is the following result which realizes the them as extensions of one another.
\begin{proposition}[{\cite[Proposition 3.3 and 3.4]{Addington}}]\label{extensions_spinors}
    Let $(V,q)$ be a quadratic space, let $W \subset V$ be an isotropic subspace and let $W' \subset W$ denote a codimension $1$ subspace. Then there exist exact sequences 
    \begin{equation}\label{possibly_interes_ext}
    \begin{aligned}
        0 \longrightarrow \mathcal{S}^W \longrightarrow \mathcal{S}^{W'} \longrightarrow \mathcal{T}^W \longrightarrow 0, \\
         0 \longrightarrow \mathcal{T}^W \longrightarrow \mathcal{T}^{W'} \longrightarrow S^W \longrightarrow 0
    \end{aligned}
    \end{equation}
    which split if and only if $W\cap K = W' \cap K$.
\end{proposition}
\begin{corollary}\label{classification_spinor_corank1}
    Let $(V,q)$ be a quadratic space and $q$ of corank $1$. Let $W \subset V$ be an isotropic subspace which is not maximal. Then the following statements hold:
    \begin{enumerate}[label=\roman*)]
        \item If $K \subset W$ and $V$ is odd dimensional, we have isomorphisms 
        \begin{equation*}
            \mathcal{S}^{W} \cong \mathcal{T}^{W} \cong (\mathcal{S}_1 \oplus \mathcal{S}_2)^{2^{d_{\mathrm{max}} - (\dim(W) +1)}}.
        \end{equation*}
        In the case $V$ is even dimensional, one replaces $\mathcal{S}_1 \oplus \mathcal{S}_2$ by $\mathcal{S}^{\oplus 2}$.
        \item If $K \not\subset W$, $\dim(W) < d_{\mathrm{max}} -1$ and $V$ is odd dimensional, then there exist isomorphisms
        \begin{equation*}
            \mathcal{S}^{W} \cong \mathcal{T}^{W} \cong (\mathcal{G}_1 \oplus \mathcal{G}_2)^{2^{d_{\mathrm{max}} - (\dim(W) +2)}};
        \end{equation*}
        here $\mathcal{G}_1$ and $\mathcal{G}_2$ denote non-trivial locally free extensions of $\mathcal{S}_1$ by $\mathcal{S}_2$ and vice versa, respectively. In the case $V$ is even dimensional, one replaces $\mathcal{G}_1 \oplus \mathcal{G}_2$ by $\mathcal{G}^{\oplus 2}$, where $\mathcal{G}$ denotes the non-trivial self extension of $\mathcal{S}$.
    \end{enumerate}
\end{corollary}
\begin{proof}
    The first case follows inductively from Proposition \ref{extensions_spinors}. For $K \not\subset W$ the claim also follows inductively from Proposition \ref{extensions_spinors}, but this inductive process doesn't start from the maximal isotropic subspace $\wmax$, but from the isotropic subspace of maximal dimension satisfying $ K \not\subset W$ which is the isotropic subspace associated with the locally free extensions $\mathcal{G}_1$ and $\mathcal{G}_2$ or $\mathcal{G}$, respectively.
\end{proof}
    The sheaves $\mathcal{G}_1, \mathcal{G}_2$ and $\mathcal{G}$ can be regarded as the ``smallest" spinor sheaves on $Y$ that are locally free and we refer to them as \textit{Kawamata's locally free extensions}. This is motivated by the case of a $5$-dimensional quadratic space $(V,q)$ with $q$ of corank $1$, where one can identify $\mathcal{G}_1, \mathcal{G}_2$ with the locally free sheaves $G_1,G_2$ appearing in \cite[Lemma 6.2]{Kawamata_SOD_ordinary_double}. Moreover, using the definitions of \cite{Kawamata_nc}, $\mathcal{G}_1 \oplus \mathcal{G}_2$ defines a \textit{$2$-pointed noncommutative deformation} of $\mathcal{S}_1\oplus\mathcal{S}_2$. 
\begin{example}
Let $(V,q)$ be a $5$-dimensional quadratic space and $q$ of corank $1$. Then, up to isomorphism, the following spinor sheaves exist on the associated quadric threefold $Y \subset \mathbb{P}(V)$:
\begin{center}
    \begin{longtable}{ | m{5cm} | m{6cm}| m{2.1cm} | } 
  \hline
  isotropic subspace & spinor sheaves & generic rank  \\ 
  \hline 
     Let $\wmax \subset V$ be the maximal isotropic subspace.
    \par
    {\centering\vspace{6pt}
	   \includegraphics[trim={15cm 34.5cm 0cm 32.3cm}, clip=true, scale = 0.18]{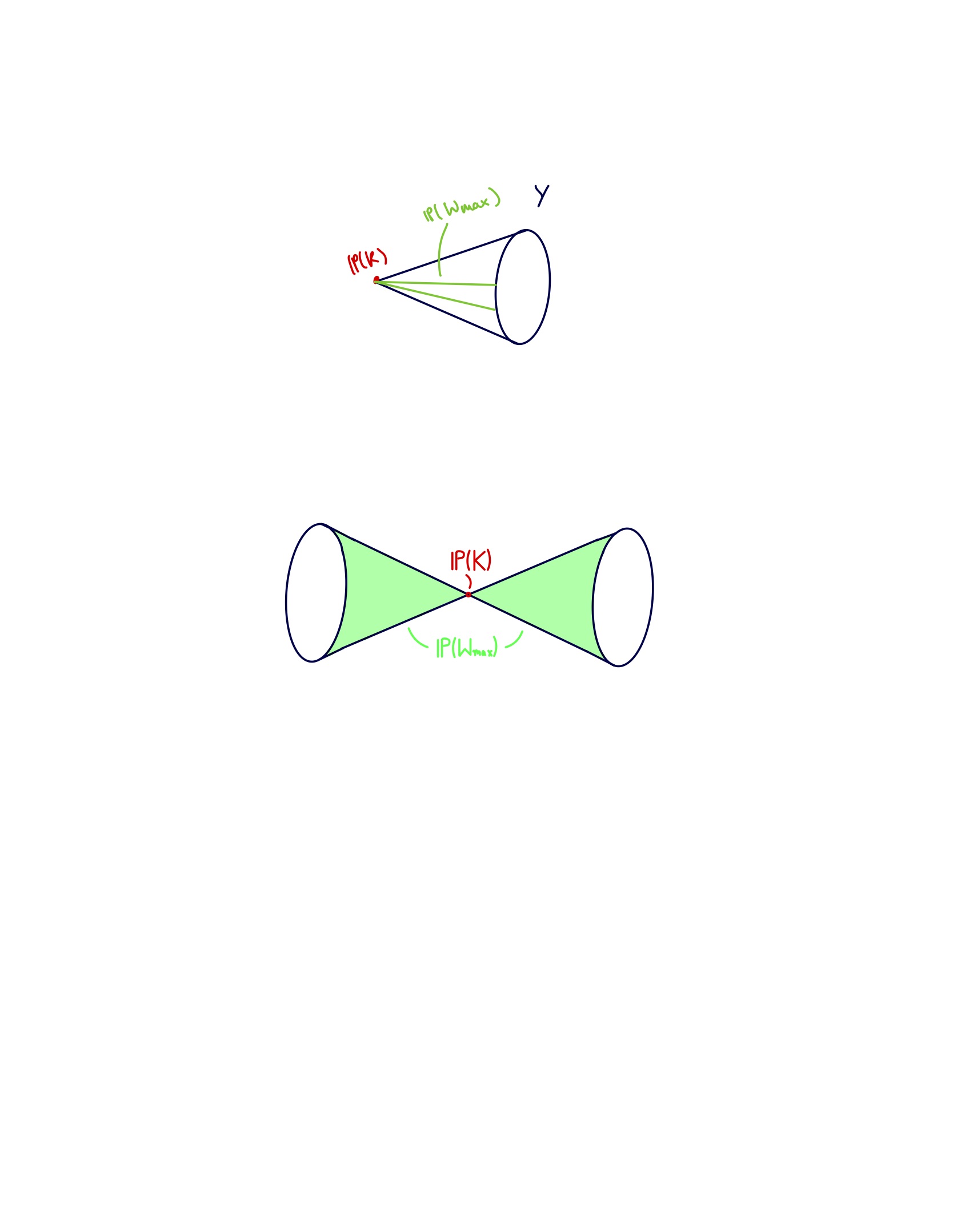}
	   \par}   & 
        We have $K \subset \wmax$. The sheaves $\mathcal{S}_1, \mathcal{S}_2$ corresponding to the maximal isotropic subspace are not locally free since they are of rank $2$ at the nodal singularity of $Y$.
        & 1\\
  \hline
  Let $W' \subset V$ be an isotropic subspace of dimension $2$. 
  \par
{\centering\vspace{6pt}
	\includegraphics[trim={12.5cm 55.5cm 0cm 9cm}, clip=true, scale = 0.18]{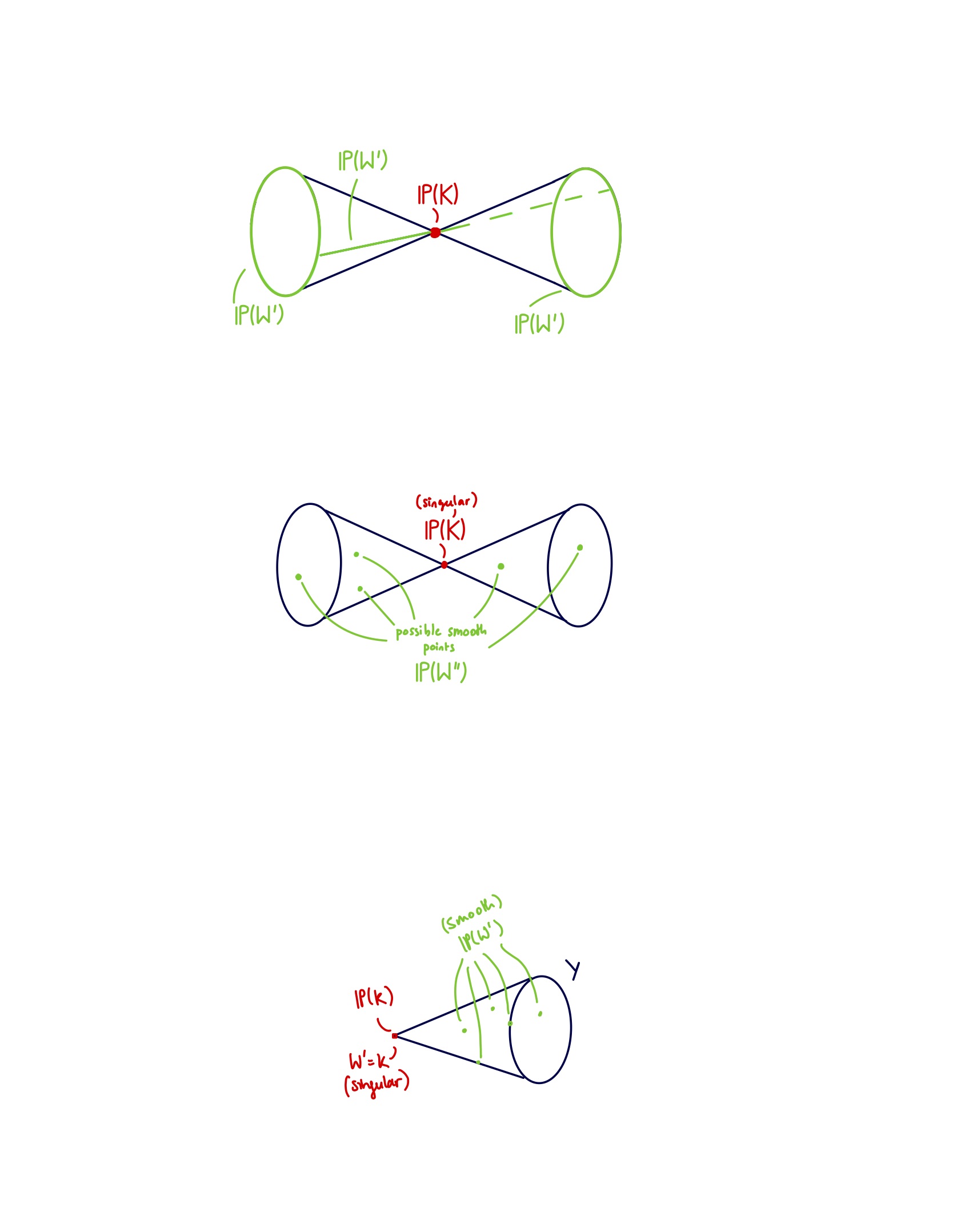}
	\par} & 
 There exist two choices for $W'$ (up to isomorphism): 
 \begin{enumerate}[label=\roman*)]
     \item If $K \subset W'$, then we have $\mathcal{S}^{W'} \cong \mathcal{T}^{W'} \cong \mathcal{S}_1 \oplus \mathcal{S}_2$;
     \item If $K \not\subset W'$, then the spinor sheaves  are Kawamata's locally free extensions $\mathcal{G}_1,\mathcal{G}_2$.
 \end{enumerate}
& 2\\
 \hline 
Let $W''\subset V$ be an isotropic subspace of dimension $1$. 
  \par
{\centering\vspace{6pt}
	\includegraphics[trim={14.5cm 34cm 0cm 30cm}, clip=true, scale = 0.19]{isotropic-subspaces.jpg}
	\par}  & 
 There exist two choices for $W''$ (up to isomorphism): 
\begin{enumerate}[label=\roman*)]
    \item If $K = W''$, then we have 
        $\mathcal{S}^{W''} \cong \mathcal{T}^{W''} \cong\mathcal{S}_1^{\oplus 2} \oplus \mathcal{S}_2^{\oplus 2}$
    \item If $K \neq W''$, then we obtain
           $\mathcal{S}^{W''} \cong \mathcal{T}^{W''} \cong \mathcal{G}_1 \oplus \mathcal{G}_2$
\end{enumerate}
& $4$ \\ 
 \hline 
    $W=\{0\}$ & \vspace{0.1cm} We have $I^W = \mathrm{Cl}(q)$. Denoting the associated spinor sheaves by 
    $\mathcal{S}_0, \mathcal{T}_0$, we have: 
    
    $\mathcal{S}_0 \cong \mathcal{T}_0 \cong \mathcal{G}_1^{\oplus 2} \oplus \mathcal{G}_2^{\oplus 2}.$
    \vspace{0.1cm}
    & $8$ \\ 
  \hline
\end{longtable}
\end{center}
\end{example}
\subsection{Clifford algebras and Morita reduction}\label{clifford_alg_Mortia}
In this subsection we will further analyze the structure of the even Clifford Algebra $\evencl(q)$, for a general nontrivial quadratic space $(V,q)$ with $q \neq 0$. If the corank of $q$ is $1$, there exist exactly one or two simple left $\evencl(q)$-modules $S_1,S_2$ and $S$, for $V$ odd or even dimensional, respectively. In this subsection the goal is to compute the corresponding Ext-complexes of Theorem \ref{heart_intro} for the simple left $\evencl(q)$-modules $S_1, S_2$ and $S$ in low dimensions.

Since we are working over an algebraically closed field $k$ with $\characteris(k) \neq 2$, the quadratic form can be written as $q= x_1^2 + \dots +x_r^2$, where $r \leq N :=\dim(V)$. For such a quadratic form we will use the notation
\begin{equation*}
    q = \langle \underbrace{1,\dots, 1}_{r\text{-times}},\underbrace{0,\dots,0}_{(N-r)\text{-times}} \rangle = \langle 1,\dots,1\rangle \perp \langle0, \dots,0\rangle,
\end{equation*}
where $\langle 1,\dots,1\rangle $ is a nondegenerate quadratic space of dimension $r$, $\langle0, \dots,0\rangle$ is a totally degenerate quadratic space of dimension $N-r$, and the symbol $\perp$ denotes their orthogonal sum. On the level of the corresponding Clifford algebra, a splitting into orthogonal sums has the following effect: Let $(V,q)$, $(V',q')$ be two arbitrary quadratic spaces. Then there exists an isomorphism of graded $k$-algebras
    \begin{equation*}
        \cl(q \perp q') \cong \cl(q) \ \widehat{\otimes} \ \cl(q'),
    \end{equation*}
see \cite[Ch. V,§1, Lemma 1.7]{Lam_quadratic}. In other words, the Clifford algebra of the orthogonal sum is isomorphic to the graded tensor product of the Clifford algebras corresponding to the respective summands. For the definition of the graded tensor product of $k$-algebras we refer to \cite[Ch. IV,§2]{Lam_quadratic}.
Analogously, one can verify a similar result for the (ungraded) even part of the Clifford algebra of an orthogonal sum of quadratic forms. 
\begin{proposition}[{\cite[Theorem 2.8]{Lam_quadratic}}]\label{Lams_theorem_odd}
    Let $(V,q)$ be an odd dimensional nondegenerate quadratic space. Then for any quadratic space $(V',q')$ there exists an isomorphism of $k$-algebras
    \begin{equation*}
        \evencl( q \perp q') \cong \evencl(q) \otimes \cl(-q').
    \end{equation*}
\end{proposition}
\begin{proof}
    Let $\{v_1,\dots,v_N\}$ denote an orthogonal basis for the odd dimensional quadratic space $(V,q)$ and set $z:=\prod^N_iv_i$. A straightforward computation shows that $z$ lies in the center of $\cl(q)$ and the claim follows from \cite[Theorem 2.8]{Lam_quadratic}.
\end{proof}
\begin{corollary}\label{lam_corr}
    Let $(V',q')$ be a quadratic space with $q\neq 0$. Then we obtain an isomorphism of (ungraded) $k$-algebras
    \begin{equation}
        \evencl(\langle 1 \rangle \perp q') \cong \cl(-q').
    \end{equation}
\end{corollary}
\begin{proof}
    The statement follows from the previous proposition. For an explicit isomorphism we refer to \cite[Lemma 4.1]{physics}.
\end{proof}
Recall that the \textit{hyperbolic plane} is a $2$-dimensional quadratic space $(V,q)$ that is isometric to $U := \langle 1,-1\rangle$. Moreover, since we work over an algebraically closed field there exists an isometry $\langle 1,-1\rangle \cong \langle 1,1\rangle$.
\begin{example}[{\cite[Example 1.5(4)]{Lam_quadratic}}]\label{clifford_alg_hyperbolic}
    Let $U$ be the hyperbolic plane and let $\widehat{M_2(k)}$ denote the $k$-algebra $M_2(k)$ equipped with the $\mathbb{Z}/2\mathbb{Z}$-grading defined by
    \begin{equation*}
        (\widehat{M_2(k)})_0 = \Big{\{} \left[ {\begin{array}{cc}
            a & 0 \\
            0 & d \\
        \end{array} } \right] \ \Big{|} \ a,d \in k \ \Big{\}}, \quad \quad
        (\widehat{M_2(k)})_1 = \Big{\{} \left[ {\begin{array}{cc}
            0 & b \\
            c & 0 \\
        \end{array} } \right] \ \Big{|} \ b,c \in k \ \Big{\}}.
    \end{equation*}
    Then there exists an isomorphism of $\mathbb{Z}/2\mathbb{Z}$-graded $k$-algebras $\cl(U) \cong \widehat{M_2(k)}$. When we consider the Clifford algebra $\cl(U)$ as an ungraded $k$-algebra, we will simply write $\matri_2(k)$.
\end{example}
\begin{proposition}\label{morita_red}
    Let $(V,q)$ be a quadratic space with $q\neq 0$ and let $U$ denote the hyperbolic plane. Then there exists an isomorphism of $k$-algebras 
    \begin{equation*}
        \evencl(q \perp U) \cong \matri_2(\evencl(q)).
    \end{equation*}
\end{proposition}
\begin{proof}
    Since $q\neq 0$ we can write $q = \langle 1 \rangle \perp q'$ and therefore Corollary \ref{lam_corr} implies that $\cl(-q') \cong \evencl(q)$. Analogously, we obtain an isomorphism of $k$-algebras 
    \begin{equation*}
        \evencl( U \perp \langle 1 \rangle) \cong \cl(-U) \cong \cl(U) \cong \matri_2(k);
    \end{equation*}
    we refer to \cite[Proposition 2.12]{Lam_quadratic} for the isomorphism in the middle. Since $U \perp \langle 1\rangle$ is an odd dimensional nondegenerate quadratic form, an application of Proposition \ref{Lams_theorem_odd} yields 
    \begin{equation*}
        \begin{aligned}
            \evencl( U \perp \langle 1 \rangle \perp q') \cong \evencl( U \perp \langle 1 \rangle) \otimes \cl(-q') &\cong \matri_2(k)  \otimes \cl(-q') \\
            &\cong \matri_2(\cl(-q')) \\
            &\cong \matri_2(\evencl(q)).
        \end{aligned}
    \end{equation*}
\end{proof}
Given a quadratic space $(V,q)$ we will denote the category of finitely generated left $\evencl(q)$-modules by $\evencl(q) \lmod$. 
\begin{corollary}\label{morita}
    Let $(V,q)$ be a quadratic space with $q\neq 0$. Then there exists an equivalence of abelian categories 
    \begin{equation}\label{morita_equiv}
        \evencl(q) \lmod \overset{\cong}{\longrightarrow} \evencl(q \perp U)\lmod.
    \end{equation}
    In other words, the $k$-algebras $\evencl(q)$ and $\evencl(q \perp U)$ are Morita equivalent.
\end{corollary}
\begin{proof}
    In the previous proposition we proved that $\evencl(q \perp U) \cong \matri_2(\evencl(q))$. Using \cite[Corollary 5.56]{rotman_hom_alg} we obtain an equivalence of abelian categories 
    \begin{equation*}
        \evencl(q) \lmod \overset{\cong}{\longrightarrow} \matri_2(\evencl(q)) \lmod, \quad M \mapsto \begin{pmatrix}M\\M\end{pmatrix}.
    \end{equation*}
\end{proof}
For the rest of this subsection the following examples of even Clifford algebras will serve as our base cases for many explicit calculations. 
\begin{example}[{\cite[Ch.\ 4, Example 1.5]{Lam_quadratic}}]\label{example_intro}
    \begin{enumerate}[label=\roman*)]
        \item Let $(V,q)$ be a $2$-dimensional quadratic space where $q$ is of corank $1$. Choose an orthogonal $k$-basis for $V$ of the form $\{v, \varepsilon \}$, where $v^2 =1$ and $\varepsilon^2 =0$. Then one can identify the Clifford algebra $\cl(q)$ with the degenerate quaternions $\big(\frac{1,0}{k}\big)$, on the generators $\{1,v,\varepsilon, j := v\varepsilon \}$, see \cite[Ch.\ IV§1]{Lam_quadratic}. In particular, the even part of the Clifford algebra is given by 
        \begin{equation*}
            \evencl(q) \cong k \oplus kj, 
        \end{equation*}
        where $j^2 = (v\varepsilon)^2 = 0$.
        \item Let $(V,q)$ be a $3$-dimensional quadratic space and let $q$ be of corank $1$. Choose an orthogonal $k$-basis $\{v_1,v_2,\varepsilon\}$ for $V$, satisfying the relations $v_1^2=1, v_2^2=1$ and $\varepsilon^2 =0$. By the notation we introduced in the beginning of this subsection, this is equivalent to writing $q= \langle 1,1,0\rangle$. Moreover, by the definition of the Clifford algebra we have $\cl(\langle 0\rangle) \cong k[\varepsilon]/\varepsilon^2$. Using previous results and examples, we obtain the following isomorphisms of $k$-algebras
        \begin{equation*}
        \begin{aligned}
            \evencl(q) & \cong (\cl(\langle 1,1\rangle) \widehat{\otimes}\cl(\langle 0 \rangle))_0 \\
            & \cong (\widehat{\matri_2(k)} \widehat{\otimes} k[\varepsilon]/\varepsilon^2)_0 \\
            & = \Big{\{} \left[ {\begin{array}{cc}
            a & b \\
            c & d \\
        \end{array} } \right] \ \Big{|} \ a,d \in k; \ b,c \in k\varepsilon \  \Big{\}} =: \left[ {\begin{array}{cc}
            k & k\varepsilon \\
            k\varepsilon & k \\
        \end{array} } \right].\\
        \end{aligned}
        \end{equation*}
    \end{enumerate}
\end{example} 
Our next goal is to determine the simple $\evencl(q)$-modules, up to isomorphism, for an arbitrary nontrivial quadratic space $(V,q)$ and prove that they coincide with the left ideals $I_0^{\wmax}, I_1^{\wmax}\subset \cl(q)$, corresponding to a maximal isotropic subspace $\wmax \subset V$, which were defined in the beginning of Subsection \ref{addington_def}. We first consider the statement for nondegenerate quadratic spaces.
\begin{proposition}\label{Artin_wed_non-deg}
    Let $(V,q)$ be a nondegenerate quadratic space and $N :=\dim(V)$. Then the following statements hold:
    \begin{enumerate}[label=\roman*)]
        \item Let $N$ be even. Then $\evencl(q)$ is a product of two isomorphic central simple algebras over $k$, in particular, we have an isomorphism of $k$-algebras 
    \begin{equation*}
        \evencl(q) \cong \matri_m(k) \times \matri_m(k),
    \end{equation*}
    where $m = 2^{\frac{N-2}{2}}$. Up to isomorphism there exist two simple left $\evencl(q)$-modules $S_1,S_2$, with $\dim_k(S_1) = \dim_k(S_2) =m$.
    \item Let $N$ be odd. Then $\evencl(q)$ is a central simple algebra and we have an isomorphism of $k$-algebras 
    \begin{equation*}
        \evencl(q) \cong \matri_m(k),
    \end{equation*}
    where $m = 2^{\frac{N-1}{2}}$. Up to isomorphism there exists one simple left $\evencl(q)$-module $S$, with $\dim_k(S) = m$.
    \end{enumerate}
Moreover, the simple modules in $\evencl(q)$ coincide with the left ideals $I_0^{\wmax}, I_1^{\wmax} \subset \cl(q)$ corresponding to the maximal isotropic subspace $\wmax \subset V$.
\end{proposition}
\begin{proof}
    By \cite[§4.2]{Lam_quadratic} we know that $\evencl(q)$ is a central simple algebra or a product of such, and the first two claims follow from an application of the Artin-Wedderburn theorem, see for example \cite[§1.3.5]{Lam_non_comm}. An elementary calculation shows that the left $\evencl(q)$-modules $I_0^{\wmax}$ and $I_1^{\wmax}$ are of dimension $m$. Suppose that $I_0^{\wmax}$ and $I_1^{\wmax}$ are not simple, then they can be decomposed into a direct sum of simple $\evencl(q)$-modules, since the module category over a semisimple ring is semisimple, 
    contradicting the fact that their dimension is $m$.
\end{proof}
Let $\Lambda$ be a finite dimensional $k$-algebra or more generally an Artin ring. Recall that the \textit{radical} $\mathfrak{r} \subset \Lambda$ is defined as the intersection of all maximal left ideals in $\Lambda$, as well as the intersection of all maximal right ideals in $\Lambda$. In particular, the radical $\mathfrak{r} \subset \Lambda$ defines a two-sided ideal.
It can be equally characterized in the following more convenient way.
\begin{lemma}[{\cite[Proposition 3.3]{auslander}}]\label{radical}
    Let $\Lambda$ be a left Artin ring and $\mathfrak{a}\subset \Lambda$ a (two-sided) ideal such that $\mathfrak{a}$ is nilpotent and $\Lambda/\mathfrak{a}$ is semisimple. Then we have $\mathfrak{a} = \mathfrak{r}$.
\end{lemma}
\begin{theorem}\label{simples_degenerate}
    Let $(V,q)$ be a quadratic space with $q\neq 0$ and let $K \subset V$ be the kernel of the quadratic form $q$. Then the simple $\evencl(q)$-modules are $I_0^{\wmax}$ and $ I_1^{\wmax}$, up to isomorphism, if $\dim(V/K)$ is even and $I_0^{\wmax} \cong I_1^{\wmax}$ if $\dim(V/K)$ is odd. 
\end{theorem}
\begin{proof}
    We choose an orthogonal $k$-basis $\{v_1, \dots, v_s, \varepsilon_1, \dots, \varepsilon_l\}$, such that 
    \begin{equation*}
        q=\langle1,\dots,1,0,\dots,0\rangle,
    \end{equation*}
    which means that the relations $v_1^2=\dots=v_s^2=1$ and $\varepsilon_1^2=\dots=\varepsilon_l^2=0$ hold. In the following $(V/K, \bar{q})$ will denote the nondegenerate quadratic space where $\bar{q} = \langle 1,\dots,1\rangle$ and we consider $\{v_1,\dots,v_s\}$ as a $k$-basis for $V/K$, that is, we have $s =\dim(V/K)$. If $s$ is even, the ideals $I_0^{\overline{W}_{\text{max}}}, I_1^{\overline{W}_{\text{max}}} $ corresponding to the maximal isotropic subspace $\overline{W}_{\text{max}} \subset V/K$ are the only simple $\evencl(\bar{q})$-modules up to isomorphism, by Proposition \ref{Artin_wed_non-deg}. To derive the analogous result for the degenerate quadratic space $(V,q)$, we consider the map 
    \begin{equation}\label{projection}
        \pi\colon\evencl(q) \twoheadrightarrow \evencl(\bar{q}), 
    \end{equation}
    which is induced by the projection $V \twoheadrightarrow V/K$, and therefore the kernel $\ker(\pi)$ is generated by the elements $\varepsilon_1, \dots,\varepsilon_l$. By Lemma \ref{radical}, there exists an isomorphism $\evencl(\bar{q}) \cong \evencl(q)/\mathfrak{r}$, since the elements in $\ker(\pi)$ are nilpotent and $\evencl(\bar{q})$ is semisimple. Therefore we have $\mathfrak{r} =\ker(\pi) = (\varepsilon_1,\dots,\varepsilon_l)$. Via the map $\pi$ we can consider the simple left $\evencl(\bar{q})$-modules $I_0^{\overline{W}_{\text{max}}}$ and $ I_1^{\overline{W}_{\text{max}}} $ as modules over $\evencl(q)$ and we will denote them by $I_{0,\evencl(q)}^{\overline{W}_{\text{max}}}$ and $I_{1,\evencl(q)}^{\overline{W}_{\text{max}}}$, respectively. They are also simple modules over $\evencl(q)$ since any $\evencl(q)$-submodule of $I_{0,\evencl(q)}^{\overline{W}_{\text{max}}}$ or $I_{1,\evencl(q)}^{\overline{W}_{\text{max}}}$ descends to a $\evencl(\bar{q})$-submodule of $I_0^{\overline{W}_{\text{max}}}$ or $ I_1^{\overline{W}_{\text{max}}}$, respectively. Moreover, the radical $\mathfrak{r}$ acts trivially on any simple left $\evencl(q)$-module M, since $\mathfrak{r}M \subset M$ is a submodule, therefore $\mathfrak{r}M = 0$ or $\mathfrak{r}M = M$, but since $\evencl(q)$ is Artinian, $\mathfrak{r}M = M$ implies $M = 0$, by Nakayama. This shows that $I_{0,\evencl(q)}^{\overline{W}_{\text{max}}}$ and $I_{1,\evencl(q)}^{\overline{W}_{\text{max}}}$ are the only simple $\evencl(q)$-modules up to isomorphism.
    
    To show that they are isomorphic to the left ideals $I_0^{\wmax}, I_1^{\wmax} \subset \cl(q)$, corresponding to the maximal isotropic subspace $\wmax \subset V$, we choose a $k$-basis of $\wmax$ of the form $\{(v_1+iv_2), (v_3+iv_4), \dots, (v_{s-1}+iv_s), \varepsilon_1,\dots,\varepsilon_l\}$. Then we have
    \begin{equation*}
        I_0^{\wmax} = \begin{cases}
            \evencl(q) (v_1+iv_2)\cdots(v_{s-1}+iv_s)\varepsilon_1\cdots\varepsilon_l, \quad \text{if } l \ \text{is even}, \\
            \oddcl(q) (v_1+iv_2)\cdots(v_{s-1}+iv_s)\varepsilon_1\cdots\varepsilon_l, \quad \text{if } l \ \text{is odd};
        \end{cases}
    \end{equation*}
    and 
    \begin{equation*}
        I_1^{\wmax} = \begin{cases}
            \oddcl(q) (v_1+iv_2)\cdots(v_{s-1}+iv_s)\varepsilon_1\cdots\varepsilon_l, \quad \text{if } l \ \text{is even}, \\
            \evencl(q) (v_1+iv_2)\cdots(v_{s-1}+iv_s)\varepsilon_1\cdots\varepsilon_l, \quad \text{if } l \ \text{is odd}.
        \end{cases}
    \end{equation*}
    Moreover, by their definition we have
    \begin{equation*}
    \begin{aligned}
        I_{0}^{\overline{W}_{\text{max}}} = \evencl(\bar{q})(v_1+iv_2)\cdots(v_{s-1}+iv_s), \\
        I_{1}^{\overline{W}_{\text{max}}} = \oddcl(\bar{q})(v_1+iv_2)\cdots(v_{s-1}+iv_s),
    \end{aligned}
    \end{equation*}
    since $\{(v_1+iv_2),\dots, (v_{s-1}+iv_s)\}$ is a basis for the maximal isotropic subspace $\overline{W}_{\text{max}}$.
    In the case where $l$ is even, we can conclude that there exist isomorphisms 
    \begin{equation}\label{isos_key}
    I_{0,\evencl(q)}^{\overline{W}_{\text{max}}}\overset{\sim}{\longrightarrow} I_0^{\wmax} \quad \text{and} \quad I_{1,\oddcl(q)}^{\overline{W}_{\text{max}}} \overset{\sim}{\longrightarrow} I_1^{\wmax},
    \end{equation}
    induced by right multiplication with the product $\varepsilon_1\cdots\varepsilon_l$. More precisely, we have a factorization
    \begin{equation*}
        \begin{tikzcd}
        \evencl(q) \arrow[r, "\varepsilon_1\cdots\varepsilon_l"] \arrow[d, dashed] & \evencl(q) \\
        \evencl(q)/\mathfrak{r} \arrow[ru, "\psi"']                                 &           
        \end{tikzcd}
    \end{equation*}
    where all the maps are surjections. The isomorphisms (\ref{isos_key}) are the restrictions of $\psi$ to $I_{0,\evencl(q)}^{\overline{W}_{\text{max}}}$ and $I_{1,\evencl(q)}^{\overline{W}_{\text{max}}}$, respectively. These maps are in fact isomorphisms because the ideals $I_0^{\wmax}, I_1^{\wmax}, I_{0,\evencl(q)}^{\overline{W}_{\text{max}}}$ and $I_{1,\evencl(q)}^{\overline{W}_{\text{max}}}$ all have the same dimension, because $\codim(\overline{W}_{\text{max}}) = \codim(\wmax)$, see Definition \ref{definition_spinor}.
    In the case where $l$ is odd, we obtain the analogous isomorphisms with $I_1^{\wmax}$ and $I_0^{\wmax}$ interchanged. If $s = \dim(V/K)$ is odd, the proof works analogously.
\end{proof}
We will state some classical results about the simple and projective modules over finite dimensional $k$-algebras, or more generally over Artin rings, from which we can construct cyclic projective resolutions of the simple $\evencl(q)$-modules.
\begin{theorem}[{\cite[§1.4.4, §1.4.5 ]{auslander}}]\label{correspondence}
    Let $\Lambda$ denote an Artin ring. Then there exists a bijection between the nonisomorphic indecomposable projective $\Lambda$-modules and the nonisomorphic simple $\Lambda$-modules. Explicitly, this is given by the map that sends an indecomposable projective $\Lambda$-module $P$ to the simple $\Lambda$-module $P/\mathfrak{r}P$.
\end{theorem}
To explicitly determine the indecomposable projective $\evencl(q)$-modules, we recall the following definitions: A set of idempotents $\{e_1,\dots,e_m\}$ is called \textit{orthogonal} if $e_ie_j = 0$ for $i \neq j$ and an idempotent element $e\in \Lambda$ is called \textit{primitive} if it cannot be written as the sum of two nonzero orthogonal idempotents. 
\begin{proposition}[{\cite[Ch.1, Proposition 4.8]{auslander}}]\label{projectives_decomp}
    Let $\Lambda$ be a left Artin ring. Then we have $1=e_1+\dots+e_m$ for suitable primitive orthogonal idempotents $e_1,\dots,e_m$ and a decomposition 
    \begin{equation*}
        \Lambda = P_1 \oplus \dots \oplus P_m,
    \end{equation*}
    where $P_i= \Lambda e_i$ are indecomposable projective $\Lambda$-modules. 
\end{proposition}
\begin{theorem}\label{extensions_and_res}
    Let $(V,q)$ be a quadratic space with $q\neq 0$ of corank $1$ and denote by $\{v_1, \dots, v_{N-1}, \varepsilon\}$ an orthogonal $k$-basis for $V$ such that $v_i^2=1$ for all $i \in \{1,\dots, N-1\}$ and $\varepsilon^2 =0$.
    \begin{enumerate}[label=\roman*)]
        \item Let $N$ be odd and let $S_1,S_2$ be the two simple $\evencl(q)$-modules which are unique up to isomorphism by Theorem \ref{simples_degenerate}. Then there exist non-split extensions of left $\evencl(q)$-modules
        \begin{equation}\label{extensions_comp}
            0 \longrightarrow S_2 \overset{\cdot\varepsilon}{\longrightarrow} P_1 \longrightarrow S_1 \longrightarrow 0
            \quad \text{and} \quad 
            0 \longrightarrow S_1 \overset{\cdot\varepsilon}{\longrightarrow} P_2 \longrightarrow S_2 \longrightarrow 0,
        \end{equation}
        where $P_1,P_2$ are the unique indecomposable projective $\evencl(q)$-modules and the exact sequences are induced by right multiplication with $\varepsilon$. In particular, there exist $2$-periodic projective resolutions 
        \begin{equation}\label{projective_resolutions}
        \begin{aligned}
            \dots \overset{\cdot\varepsilon}{\longrightarrow} P_2 \overset{\cdot\varepsilon}{\longrightarrow} P_1 \overset{\cdot\varepsilon}{\longrightarrow}P_2 \overset{\cdot\varepsilon}{\longrightarrow} P_1 \longrightarrow S_1;\\
            \dots \overset{\cdot\varepsilon}{\longrightarrow} P_1\overset{\cdot\varepsilon}{\longrightarrow} P_2 \overset{\cdot\varepsilon}{\longrightarrow} P_1 \overset{\cdot\varepsilon}{\longrightarrow} P_2 \longrightarrow S_2.
        \end{aligned}
        \end{equation} 
        \item Let $N$ be even and let $S$ denote the simple $\evencl(q)$-module which is unique up to isomorphism by Theorem \ref{simples_degenerate}. Set $j := v_{N-1}\varepsilon$.
        Then there exists a non-split extension of left $\evencl(q)$-modules 
        \begin{equation}
            0 \longrightarrow S \overset{\cdot j}{\longrightarrow} P \longrightarrow S\longrightarrow 0,
        \end{equation}
        where $P := \evencl(q)$ is the unique indecomposable projective $\evencl(q)$-module, up to isomorphism and the exact sequence is induced by right multiplication with $j$. In particular, there exists a $1$-periodic projective resolution
        \begin{equation}\label{even_proj_res}
            \dots \overset{\cdot j}{\longrightarrow} P \overset{\cdot j}{\longrightarrow} P \overset{\cdot j}{\longrightarrow} P \overset{\cdot j}{\longrightarrow} P \longrightarrow S.
        \end{equation}
        \end{enumerate}
\end{theorem}
\begin{proof}
    Since simple modules are preserved under the Morita equivalence (\ref{morita_equiv}), it suffices to prove the statement for the low dimensional forms $q= \langle 1,1,0\rangle$ and $q=\langle 1,0 \rangle$ for $N$ odd and even, respectively. 
    
    Starting with the odd case, consider Example \ref{example_intro}(ii), that is, there exists an isomorphism 
    \begin{equation*}
        \evencl(q) \cong \left[ {\begin{array}{cc}
            k & k\varepsilon \\
            k\varepsilon & k \\
        \end{array} } \right].
    \end{equation*}
    With Lemma \ref{radical} we can deduce that the radical of $\evencl(q)$ is given by 
    \begin{equation*}
        \mathfrak{r} = \left[ {\begin{array}{cc}
            0 & k\varepsilon \\
            k\varepsilon & 0 \\
        \end{array} } \right].
    \end{equation*}
    To determine the indecomposable projective $\evencl(q)$-modules we can pick the primitive orthogonal idempotents 
    \begin{equation*}
        e_1 = \left[ {\begin{array}{cc}
            1 & 0 \\
            0 & 0 \\
        \end{array} } \right], \quad 
        e_2 = \left[ {\begin{array}{cc}
            0 & 0 \\
            0 & 1 \\
        \end{array} } \right] 
    \end{equation*}
    which satisfy $1 = e_1 + e_2$. By Proposition \ref{projectives_decomp}, the indecomposable projective left $\evencl(q)$-modules are given by the column vectors 
    \begin{equation*}
        P_1 = \evencl(q)e_1 = \begin{pmatrix}k\\k\varepsilon\end{pmatrix},  \quad 
        P_2 = \evencl(q)e_2 = \begin{pmatrix}k\varepsilon\\k\end{pmatrix}
    \end{equation*}
    and together with Theorem \ref{correspondence} we obtain that the simples are given by 
    \begin{equation*}
        S_1 = \begin{pmatrix}k\\0\end{pmatrix}, \quad S_2 = \begin{pmatrix}0\\k\end{pmatrix}.
    \end{equation*}
    This gives rise to nontrivial extensions of $\evencl(q)$-modules 
    \begin{equation*}
        0 \longrightarrow S_2 \overset{\cdot\varepsilon}{\longrightarrow} P_1 \longrightarrow S_1 \longrightarrow 0
        \quad \text{and} \quad 
        0 \longrightarrow S_1 \overset{\cdot\varepsilon}{\longrightarrow} P_2 \longrightarrow S_2 \longrightarrow 0,
    \end{equation*}
    induced by right multiplication with $\varepsilon$.
    Splicing these sequences together yields the desired projective resolutions of $S_1$ and $S_2$. 
    
    For $N$ even, we consider the quadratic form $q = \langle 1,0\rangle$, as we did in Example \ref{example_intro}(i), from which we know that $\evencl(q) = k \oplus kj$, where $j^2=0$. 
    An application of Lemma \ref{radical} shows that the 
    ideal $kj \subset \evencl(q)$ coincides with the radical ideal. By the same arguments as in the odd dimensional case we can conclude that the only simple module is given by $S =k \subset \evencl(q)$ and the only indecomposable projective $\evencl(q)$-module is given by $\evencl(q)$ itself. We obtain a non-split extension
    \begin{equation*}
        0 \longrightarrow k \overset{\cdot j}{\longrightarrow} \evencl(q) \longrightarrow k\longrightarrow 0
    \end{equation*}
    which yields the desired projective resolution of $S$.
\end{proof}
\begin{proposition}\label{technical_heart}
   In the setting of Theorem \ref{extensions_and_res} the following statements hold.
    \begin{enumerate}[label=\roman*)]
        \item Let $N$ be odd dimensional. Then there exist isomorphisms of $k$-algebras
        \begin{equation*}
        \begin{aligned}
            \mathrm{Ext}^{\bullet}(S_1,S_1) \cong \mathrm{Ext}^{\bullet}(S_2,S_2)  \cong k[\theta],
        \end{aligned}
    \end{equation*}
    where $\theta$ denotes an element of degree $2$.
    \item Let $N$ be even dimensional. Then there exists an isomorphism of $k$-algebras
        \begin{equation*}
        \begin{aligned}
            \mathrm{Ext}^{\bullet}(S,S) \cong k[\theta'],
        \end{aligned}
    \end{equation*}
    where $\theta'$ denotes an element of degree $1$.
    \end{enumerate}
\end{proposition}
\begin{proof}
    Let $N$ be odd dimensional. Considering Theorem \ref{extensions_and_res}(i) we have 
    \begin{equation*}
        \Hom(P_2, S_1) \cong 0 \cong \Hom(P_1, S_2),
    \end{equation*}
    since the extensions (\ref{extensions_comp}) are nontrivial. Moreover, we have isomorphisms
    \begin{equation*}
         \Hom(P_1, S_1) \cong k \cong \Hom(P_2, S_2),
    \end{equation*}
    that follow from the fact that $S_1, S_2$ are simple. More precisely, since $\Hom(S_2,S_1) =0$ and $S_2 \overset{\cdot \varepsilon}{\hookrightarrow} P_1$, we see that $\Hom(P_1, S_1) \cong \Hom(S_1, S_1) \cong k$, because every morphism $P_1 \to S_1$ uniquely factors through a morphism $P_1/S_2\cong S_1 \to S_1$.
    
    We apply the functor $\Hom(-, S_1)$ and $\Hom(-,S_2)$ to the projective resolutions (\ref{projective_resolutions}), respectively. Then, taking  cohomology gives rise to isomorphisms of graded $k$-vector spaces 
    \begin{equation*}
        \mathrm{Ext}^{\bullet}(S_1, S_1) \cong k[\theta] \cong \mathrm{Ext}^{\bullet}(S_2, S_2),
    \end{equation*}
    where $\theta$ denotes a generator of degree $2$.
    To prove that they are in fact isomorphisms of $k$-algebras, we consider the exact triangles
    \begin{equation}\label{sequences_algebra}
        \begin{aligned}
        S_2 \longrightarrow P_1 \longrightarrow S_1 \overset{\eta}{\longrightarrow} S_2[1]
        \quad \text{and} \quad 
        S_1 \longrightarrow P_2 \longrightarrow S_2 \overset{\eta'}{\longrightarrow} S_1[1],
        \end{aligned}
    \end{equation}
    where $\eta, \eta'$ denote nontrivial morphisms of degree $1$ corresponding to the nontrivial extensions $P_1$ and $P_2$. This gives rise to triangles of graded $k[\theta]$-modules
    \begin{equation*}
    \begin{aligned}
        \mathrm{Ext}^{\bullet}(S_2,S_1)[-1] \overset{\eta}{\longrightarrow} \mathrm{Ext}^{\bullet}(S_1,S_1) \longrightarrow k; \\
        \mathrm{Ext}^{\bullet}(S_1,S_1)[-1] \overset{\eta'}{\longrightarrow} \mathrm{Ext}^{\bullet}(S_2,S_1) \longrightarrow k.
    \end{aligned}
    \end{equation*}
    In particular, the maps induced by $\eta$ and $\eta'$ are injective. Therefore, composition with the element $\eta'[1]\circ\eta \in \mathrm{Ext}^{2}(S_1,S_1)$ yields an injective morphism 
    \begin{equation*}
        \mathrm{Ext}^{\bullet}(S_1,S_1)[-2] \longrightarrow \mathrm{Ext}^{\bullet}(S_1,S_1).
    \end{equation*}
    We obtain an isomorphism of $k$-algebras
    \begin{equation*}
        k[\theta] \overset{\cong}{\longrightarrow} \mathrm{Ext}^{\bullet}(S_1,S_1), \quad \theta \mapsto (\eta'[1] \circ\eta\colon S_1 \to S_1[2]).
    \end{equation*}
    To compute $\mathrm{Ext}^{\bullet}(S_2,S_2)$, one applies $\Hom(-,S_2)$ to the exact triangles (\ref{sequences_algebra}) and the statement follows by the same arguments. \\
    Let $N$ be even dimensional and consider the projective resolution (\ref{even_proj_res}). 
    Applying the functor $\Hom(-,S)$ and taking cohomology yields an isomorphism of $k$-vector spaces
    \begin{equation*}
        \mathrm{Ext}^{\bullet}(S,S) \cong k[\theta'], 
    \end{equation*}
    where $\theta'$ denotes a generator of degree $1$. To see that this is in fact an isomorphism of $k$-algebras, we consider the exact triangle 
    \begin{equation*}
         S \longrightarrow \evencl(q) \longrightarrow S\overset{\nu}{\longrightarrow} S[1],
    \end{equation*}
    where $0\neq \nu \in \mathrm{Ext}^1(S,S)$ corresponds to the nontrivial extension $P$. This gives rise to an exact triangle
    \begin{equation*}
        \mathrm{Ext}^{\bullet}(S,S)[-1] \overset{\nu}{\longrightarrow} \mathrm{Ext}^{\bullet}(S,S) \longrightarrow k
    \end{equation*}
    which induces an isomorphism of $k$-algebras 
    \begin{equation*}
        k[\theta'] \overset{\cong}{\longrightarrow}\mathrm{Ext}^{\bullet}(S,S), \quad \theta' \mapsto (\nu\colon S \to S[1]).
    \end{equation*}
\end{proof}
Let $(V,q)$ be an odd dimensional quadratic space and $q\neq 0$ of corank $1$. For any module $M \in D^b(\evencl(q))$, the complex $\mathrm{Ext}^{\bullet}(M,S_1)$ admits a $k[\theta]$-action by postcomposing with the generator $\theta\colon S_1 \to S_1[2]$. For the following statement, let us consider the non-split extension 
\begin{equation}\label{ses}
        0 \longrightarrow S_1 \longrightarrow P_2 \longrightarrow S_2 \longrightarrow 0,
\end{equation}
of Theorem \ref{extensions_and_res} as an element $0\neq \kappa \in \Hom(S_2, S_1[1]) \cong \mathrm{Ext}^{1}(S_2, S_1)$.
\begin{lemma}\label{module-structure}
    Let $(V,q)$ be an odd dimensional quadratic space and $q\neq 0$ a quadratic form of corank $1$. Then the complex $\mathrm{Ext}^{\bullet}(S_2, S_1)$ admits a structure of a free $k[\theta]$-module of rank $1$, generated by the element $\kappa \in \mathrm{Ext}^{1}(S_2, S_1)$.
\end{lemma}
\begin{proof}
    Since $P_2$ is projective, applying $\Hom(-,\mathcal{S}_1)$ to the short exact sequence (\ref{ses})
    induces an isomorphism of $k[\theta]$-modules
	\begin{equation}\label{k_theta_module}
		\kappa^*\colon\mathrm{Ext}^{\bullet}(S_1, S_1)[-1] \overset{\sim}{\longrightarrow}\mathrm{Ext}^{\bullet}(S_2,S_1), \quad f \mapsto f\circ \kappa.
	\end{equation}
\end{proof}
\begin{remark}\label{k_theta_k_theta_bimodule_str}
    Note that the complex $\mathrm{Ext}^{\bullet}(S_2,S_1)$ admits the structure of a $k[\theta]$-$k[\theta]$-bimodule, by precomposing with the generator  $\theta\colon S_2[-2] \to S_2$. 
\end{remark}
\begin{remark}\label{module_str_other}
    We can prove the analogous result for $\mathrm{Ext}^{\bullet}(S_1, S_2)$, by using the extension 
    \begin{equation}
        0 \longrightarrow S_2 \longrightarrow P_1 \longrightarrow S_1 \longrightarrow 0,
    \end{equation}
    established in Theorem \ref{extensions_and_res}.
\end{remark}
\subsection{The derived category of a quadric and spinor sheaves}\label{derived_spinor}
Let $(V,q)$ be a quadratic space with $q$ of corank $1$. In this subsection we will establish an equivalence between the derived category of left $\evencl(q)$-modules and the triangulated subcategory of $D^b(Y)$ generated by the sheaves $\mathcal{S}_1,\mathcal{S}_2$ and $\mathcal{S}$ for an odd and even dimensional nodal quadric $Y := V(q) \subset \mathbb{P}(V)$, respectively. 

In the following, $\mathcal{S}_0$ and $\mathcal{T}_0$ will denote the spinor sheaves corresponding to the trivial isotropic subspace $\{0\}\subset V$. These can be seen as the ``biggest" spinor sheaves with respect to their rank, and they are always locally free, see Proposition \ref{locallyfree_or_not}. In addition, the spinor sheaves $\mathcal{S}_0,\mathcal{T}_0$ carry a natural right $\evencl(q)$-module structure, since the morphisms defining $\mathcal{S}_0$ and $\mathcal{T}_0$, are induced by left multiplication
    \begin{equation*}
        \evencl(q) \overset{v\cdot}{\longrightarrow} \oddcl(q) \quad \text{and} \quad \oddcl(q) \overset{v\cdot}{\longrightarrow} \evencl(q),
    \end{equation*}
for a vector $v\in V$.

Now we can define the functor inducing the equivalence we described above: 
\begin{proposition}[{\cite[Theorem 4.2]{quadric_fibr}}]\label{Kuznetsov_ff_SOD}
    Let $(V,q)$ be a quadratic space, let $Y \subset \mathbb{P}(V)$ be the associated quadric hypersurface of dimension $n$ and let $p\colon \mathbb{P}(V) \to k$ be the natural projection. Then there exists a fully faithful functor
    \begin{equation}
        \Phi\colon D^b(\evencl(q)) \longrightarrow D^b(Y), \quad M \mapsto \mathcal{S}_0 \otimes_{\evencl(q)} p^*M.
    \end{equation}
    Moreover, there exists a semiorthogonal decomposition
    \begin{equation}
        D^b(Y) = \langle D^b(\evencl(q)), \mathcal{O},\dots,  \mathcal{O}(n-1) \rangle
    \end{equation}
    which yields an equivalence between the Kuznetsov component $\mathcal{A}_Y$ and the category $D^b(\evencl(q))$. 
\end{proposition}
\begin{remark}
    In \cite{quadric_fibr} the author assumes $D^b(\evencl(q))$ to be the derived category of \textit{right} $\evencl(q)$-modules. Since we work in the setting of left $\evencl(q)$-modules, we rephrased his results accordingly. Moreover, we note that the sheaf $\mathcal{E}'_{-1,1}$, defined in \cite[preceding Lemma 4.5]{quadric_fibr}, coincides with the spinor sheaf $\mathcal{S}_0$ associated to the trivial isotropic subspace $\{0\}\subset V$, see \cite[Lemma 4.7]{quadric_fibr}. 
\end{remark}
In the following we show that the functor $\Phi$ already exists on the level of abelian categories. For a quadratic space $(V,q)$ we denote the abelian category of $\mathbb{Z}/2\mathbb{Z}$-graded left $\cl(q)$-modules by $\mathrm{gr}_{\mathbb{Z}/2\mathbb{Z}}\cl(q)\lmod$.
\begin{theorem}[{\cite[§3]{Addington}}]\label{Addingtons_functor}
    Let $(V,q)$ be a quadratic space and let $Y := V(q) \subset \mathbb{P}(V)$ denote the associated quadric hypersurface. The functor 
    \begin{equation*}
    \Psi\colon  \mathrm{gr}_{\mathbb{Z}/2\mathbb{Z}}\cl(q)\lmod  \longrightarrow \coh(Y)
    \end{equation*}
    given by 
    \begin{equation*}
    M=M_0\oplus M_1  \mapsto \coker(\mathcal{O}_{\mathbb{P}(V)}(-2) \otimes M_0 \overset{\varphi}{\to} \mathcal{O}_{\mathbb{P}(V)}(-1) \otimes M_1)
    \end{equation*}
    defines a fully faithful embedding.
\end{theorem}
The two functors $\Psi$ and $\Phi$ are related via Morita equivalence, precisely we have:
\begin{lemma}\label{bridge_Add_Kuz}
    Let $(V,q)$ be a quadratic space, such that $q\neq 0$. There exists an equivalence of categories
    \begin{equation*}
        \mathrm{gr}_{\mathbb{Z}/2\mathbb{Z}}\cl(q)\lmod  \overset{\cong}{\longrightarrow}\evencl(q)\lmod, \quad M \mapsto M_0.
    \end{equation*}
\end{lemma}
\begin{proof}
    One can show that the mapping rule $N \mapsto N \otimes_{\evencl(q)} \cl(q)$ defines a quasi-inverse of the above functor. For an explicit verification we refer to \cite[Lemma 4.3]{physics}.
\end{proof}
We can summarize the previous results in the following diagram of triangulated categories.
\begin{equation*}
\begin{tikzcd}[row sep = 3cm, column sep = 4cm]
D^b(Y) &   \\     
D^b(\mathrm{gr}_{\mathbb{Z}/2\mathbb{Z}}\cl(q)) \arrow[u, hook, "\Psi"] \arrow[r, "\text{Lemma \ref{bridge_Add_Kuz}}", "\cong"']\arrow[d,"\cong"'] & D^b(\mathrm{Cl}_0(q)) \arrow[d,"\text{Corollary \ref{morita}}", "\cong"'] \arrow[l] \arrow[lu, hook',  "\Phi:=  \mathcal{S}_0 \otimes_{\mathrm{Cl}_0(q)} p^*(-)"', bend right=15] \\
D^b(\mathrm{gr}_{\mathbb{Z}/2\mathbb{Z}}\cl(q\perp U)) \arrow[r, "\text{Lemma \ref{bridge_Add_Kuz}}", "\cong"'] \arrow[u] & D^b(\mathrm{Cl}_0(q \perp U)). \arrow[l] \arrow[u]
\end{tikzcd}
\end{equation*}
We observe that the upper triangle commutes, since we can rewrite $\Phi$ as a functor defined by the assignment 
\begin{equation*}
    N \longmapsto \cone( \mathcal{O}(-2) \otimes\mathrm{Cl}_0(q) \overset{\varphi}{\longrightarrow} \mathcal{O}(-1) \otimes \mathrm{Cl}_1(q)) \otimes_{\mathrm{Cl}_0(q)} p^*N,
\end{equation*}
using the definition of $\mathcal{S}_0$. Therefore we conclude that the functor $\Phi$ already exists on the level of abelian categories.
\begin{proposition}\label{Mod_sheaf_equivalence}
    Let $(V,q)$ be a quadratic space and $q\neq 0$ a quadratic form of corank $1$. If the associated nodal quadric $Y \subset \mathbb{P}(V)$ is odd dimensional, there exists an equivalence of categories 
    \begin{equation*}
        \Phi\colon D^b(\evencl(q)) \overset{\cong}{\longrightarrow} \langle \mathcal{S}_1, \mathcal{S}_2 \rangle \subset D^b(Y).
    \end{equation*}
    In the case $Y$ is even dimensional, we obtain an equivalence of categories
    \begin{equation*}
        \Phi\colon D^b(\evencl(q)) \overset{\cong}{\longrightarrow} \langle \mathcal{S} \rangle \subset D^b(Y).
    \end{equation*}
\end{proposition}
\begin{proof}
    The functor $\Phi$ is fully faithful, and from its definition we immediately deduce that the ideals $I_0^{\wmax}$ and $ I_1^{\wmax}$ are mapped to the spinor sheaves $\mathcal{S}_1$ and $\mathcal{S}_2$, respectively. This implies that the triangulated subcategory $\langle \mathcal{S}_1, \mathcal{S}_2 \rangle \subset D^b(Y)$ is contained in the image of $\Phi$. For the other inclusion, we proved in Theorem \ref{simples_degenerate} that the ideals $I_0^{\wmax}, I_1^{\wmax}$ are the unique simple $\evencl(q)$-modules, up to isomorphism, which implies that they generate the triangulated category $D^b(\evencl(q))$.
\end{proof}
\begin{proof}[Proof of Theorem \ref{heart_intro}]
The claims follow from Proposition \ref{technical_heart}, Lemma \ref{module-structure}, Remark \ref{module_str_other} and the equivalence of the previous proposition.
\end{proof}
With Theorem \ref{heart_intro} at hand, we can prove the following lemmas 
which will be relevant for the proof of Theorem \ref{MAIN}. Let $Y \subset \mathbb{P}(V)$ be an odd dimensional nodal quadric. Then for any sheaf $\mathcal{F} \in \coh(Y)$, the complex $\mathrm{Ext}^{\bullet}(\mathcal{F},\mathcal{S}_1)$ admits a $k[\theta]$-action by postcomposing with the generator $\theta\colon \mathcal{S}_1 \to \mathcal{S}_1[2]$. In Lemma \ref{exactsequences_spinors_o} we introduced an extension
	\begin{equation}\label{ses_alpha}
		0 \longrightarrow \mathcal{S}_2 \longrightarrow \mathcal{O}_Y^{\oplus M} \longrightarrow \mathcal{S}_1(1)\longrightarrow 0 
	\end{equation}
    which we will consider as an element $0\neq \alpha \in \Hom(\mathcal{S}_1(1), \mathcal{S}_2[1]) \cong \mathrm{Ext}^{1}(\mathcal{S}_1(1), \mathcal{S}_2)$. 
    In Lemma \ref{module-structure} the element $\kappa \in \mathrm{Ext}^1(S_2, S_1)$ represents the non-split extension of left $\evencl(q)$-modules
    \begin{equation}
        0 \longrightarrow S_1 \longrightarrow P_2 \longrightarrow S_2 \longrightarrow 0
    \end{equation}
    which is mapped to
    \begin{equation}\label{extension_kap}
        0 \longrightarrow \mathcal{S}_1 \longrightarrow \mathcal{G}_2 \longrightarrow \mathcal{S}_2 \longrightarrow 0
    \end{equation}
    under the equivalence $\Phi$ of Proposition \ref{Mod_sheaf_equivalence}. Here, $\mathcal{G}_1, \mathcal{G}_2$ denote Kawamata's locally free extensions, see Subsection \ref{addington_def}.
    In the following, we denote the class in $\mathrm{Ext}^1(\mathcal{S}_2,\mathcal{S}_1)$ corresponding to the extension (\ref{extension_kap}) also by $\kappa$.
\begin{lemma}\label{module_structure_twist}
    Let $(V,q)$ be an odd dimensional quadratic space and $q\neq 0$ a quadratic form of corank $1$. Then the complex $\mathrm{Ext}^{\bullet}(\mathcal{S}_1(1), \mathcal{S}_1)$ admits a structure of a free $k[\theta]$-module of rank $1$, generated by the element $\kappa[1] \circ \alpha \in \mathrm{Ext}^2(\mathcal{S}_1(1), \mathcal{S}_1)$.
\end{lemma}
\begin{proof}
    By Lemma \ref{cohomology}, the cohomology groups $H^i(Y, \mathcal{S}_1)$ vanish for all $i \geq 0$ and therefore the short exact sequence (\ref{ses_alpha}) gives rise to an isomorphism of $k[\theta]$-modules 
	\begin{equation*}
		\alpha^*\colon \mathrm{Ext}^{\bullet}(\mathcal{S}_2, \mathcal{S}_1)[-1]\overset{\sim}{\longrightarrow}\mathrm{Ext}^{\bullet}(\mathcal{S}_1(1), \mathcal{S}_1), \quad f \mapsto f\circ \alpha.
	\end{equation*}
    Together with the results of Lemma \ref{module-structure} this gives rise to an isomorphism of $k[\theta]$-modules
    \begin{equation}\label{k_theta_module_twist}
        \mathrm{Ext}^{\bullet}(\mathcal{S}_1, \mathcal{S}_1)[-2]\overset{\sim}{\longrightarrow}\mathrm{Ext}^{\bullet}(\mathcal{S}_1(1), \mathcal{S}_1), \quad g \longmapsto g \circ \kappa[1] \circ \alpha.
    \end{equation}
\end{proof}
If $Y \subset \mathbb{P}(V)$ is an even dimensional nodal quadric, the complex $\mathrm{Ext}^{\bullet}(\mathcal{S}(1), \mathcal{S})$ admits a $k[\theta']$-action by postcomposing with the generator $\theta'\colon \mathcal{S} \to \mathcal{S}[1]$. In Lemma \ref{exactsequences_spinors_o} we proved that there exists an extension
    \begin{equation*}
		0 \longrightarrow \mathcal{S} \longrightarrow \mathcal{O}_Y^{\oplus M} \longrightarrow \mathcal{S}(1)\longrightarrow 0
	\end{equation*}
 which we will consider as an element $0 \neq \beta \in \Hom(\mathcal{S}(1), \mathcal{S}[1]) \cong \mathrm{Ext}^{1}(\mathcal{S}(1), \mathcal{S})$.
\begin{lemma}\label{even_module_structure_twist}
    Let $(V,q)$ be an even dimensional quadratic space and $q\neq 0$ a quadratic form of corank $1$. Then the complex $\mathrm{Ext}^{\bullet}(\mathcal{S}(1), \mathcal{S})$ admits a structure of a free $k[\theta']$-module of rank $1$, generated by the element $\beta \in \mathrm{Ext}^{1}(\mathcal{S}(1), \mathcal{S})$.
\end{lemma}
\begin{proof}
    Since the cohomology groups $H^i(Y, \mathcal{S})$ vanish for all $i \geq 0$ by Lemma \ref{cohomology}, we obtain an isomorphism of $k[\theta']$-modules 
    \begin{equation}\label{even_module_str_iso}
		\beta^*\colon \mathrm{Ext}^{\bullet}(\mathcal{S}, \mathcal{S})[-1]\overset{\sim}{\longrightarrow}\mathrm{Ext}^{\bullet}(\mathcal{S}(1), \mathcal{S}), \quad f \mapsto f\circ \beta.
	\end{equation} 
\end{proof}
\subsection{An alternative proof of Theorem \ref{heart_intro} (after the anonymous referee)}\label{proof_referee}
Recall the following alternative definition of spinor sheaves on a nodal quadric which is fundamental to their proof. Let $(V,q)$ be a quadratic space with $q$ of corank $1$, denote by $Y\subset \mathbb{P}(V)$ the associated nodal quadric and the kernel of the quadratic form $q\colon V \to k$ by $K$. Then we can consider a diagram 
\begin{equation*}\label{setting}
	\begin{tikzcd}
                & E \arrow[r] \arrow[ld]         & \widetilde{Y} = \bl_y(Y) \arrow[ld, "\pi"'] \arrow[rd, "\rho"] &      \\
\{y\} \arrow[r] & Y \arrow[rr, "\rho_0", dashed] &                                                                & {Q}
\end{tikzcd}
\end{equation*}
which is induced by projecting away from the nodal point $y \in Y$ onto a smooth quadric $Q \subset \mathbb{P}(V/K)$ and extending this rational map via the blow-up $\widetilde{Y}$ of $Y$ at the nodal point $y\in Y$. We denote by $E$ the exceptional divisor of the blow-up $\pi$ and we have $\dim(Q) = \dim(Y) -1$. 
There exist spinor sheaves on the smooth quadric $Q$ corresponding to the maximal isotropic subspace in $V/K$ which we will denote by $\bar{\mathcal{S}}_1,\bar{\mathcal{S}}_2$ and $\bar{\mathcal{S}}$, for $\dim(V/K)$ even and odd, respectively. Note that by an application of Lemma \ref{cohomology} and \cite[Theorem 3.5]{Ottaviani_later} one can identify these bundles with Ottaviani's spinor bundles introduced in \cite[Definition 1.3]{Ottaviani_spinors}. 

This lets us define the spinor sheaves on Y as pullbacks of spinor bundles on the smooth quadric $Q$. More precisely, following \cite[Proposition 6.4]{categorical_absorptions} or \cite[Remark 6.6]{Kawamata_SOD_ordinary_double}, the complexes 
\begin{equation*}
    \mathcal{S}_1 := \pi_*\rho^*\bar{\mathcal{S}}_1, \quad \mathcal{S}_2 := \pi_*\rho^*\bar{\mathcal{S}}_2, \quad \text{and} \quad \mathcal{S} := \pi_*\rho^*\bar{\mathcal{S}}
\end{equation*}
are maximal Cohen--Macaulay sheaves on the odd and even dimensional nodal quadric $Y$, respectively. As the abuse of notation suggests, these coincide with the spinor sheaves on $Y$, introduced in Definition \ref{defi_spinors_sec3}, see \cite[Proposition 3.8]{A2_v2}.

The referee proposed the following alternative strategy for proving Theorem \ref{heart_intro}. We outline it in the case of an odd dimensional nodal quadric $Y$ (the even case works analogously). There exist short exact sequences 
    \begin{equation*}
        0 \to \rho^*\bar{\mathcal{S}}_2 \to \pi^*\mathcal{G}_1 \to \rho^*\bar{\mathcal{S}}_1 \otimes \mathcal{O}(-E) \to 0 \ \ \text{and} 
        \ \ 0 \to \rho^*\bar{\mathcal{S}}_1 \to \pi^*\mathcal{G}_2 \to \rho^*\bar{\mathcal{S}}_2 \otimes \mathcal{O}(-E) \to 0, 
    \end{equation*}
    where $\mathcal{G}_1$ and $\mathcal{G}_2$ denote Kawamata's locally free extensions. They are the key to determining the complexes 
    \begin{equation*}
        \mathrm{Ext}^{\bullet}(\mathcal{G}_i, \mathcal{S}_j) \cong \mathrm{Ext}^{\bullet}(\mathcal{G}_i, \pi_*\rho^*\bar{\mathcal{S}}_j) = \mathrm{Ext}^{\bullet}(\pi^*\mathcal{G}_i, \rho^*\bar{\mathcal{S}}_j) \cong \begin{cases}
        k, \quad \text{if } i =j, \\
        0, \quad \text{otherwise.}
        \end{cases}
    \end{equation*}
    By splicing together the short exact sequences of Proposition \ref{extensions_spinors} we obtain periodic long exact sequences 
    \begin{equation*}
    \begin{aligned}
        \dots \to \mathcal{G}_2 \to \mathcal{G}_1 \to \mathcal{G}_2 \to \mathcal{G}_1 \to \mathcal{S}_1 \to 0; \\
        \dots \to \mathcal{G}_1 \to \mathcal{G}_2 \to \mathcal{G}_1 \to \mathcal{G}_2 \to \mathcal{S}_2 \to 0
    \end{aligned}
    \end{equation*}
    which can be used to obtain the same Ext-complexes as in Theorem \ref{heart_intro}.  
\section{Kernels of crepant categorical resolutions of cuspidal singularities}\label{constr_crr}
The goal of this section is to prove our main theorem (Theorem \ref{MAIN}). The following result is a generalization of \cite[Proposition 3.5]{Kernels_nodal} and \cite[Theorem 5.8]{categorical_absorptions}, who prove the statement for $A_1$ singularities. 
\begin{theorem}\label{cuspcatres}
	Let $X$ be a projective variety with an isolated $A_1$ or $A_2$ singularity and $\dim(X) \geq 2$. Then there exists a crepant categorical resolution $\widetilde{\mathcal{D}}$ of $\dbx$. More precisely, we have 
    \begin{equation*}
        \widetilde{\mathcal{D}} = \{ F \in D^b(\widetilde{X}) | \ j^*F \in \langle \mathcal{A}_Y, \mathcal{O}_Y \rangle \},
    \end{equation*}
    where $\mathcal{A}_Y$ denotes the Kuznetsov component of the exceptional divisor $Y$ of the blow-up $\widetilde{X}$ of $X$ at the singularity $x$.
\end{theorem}
\begin{remark}
    In the case that $X$ is odd dimensional with an isolated $A_1$ singularity, the crepant categorical resolution $\widetilde{\mathcal{D}}$ above is ``bigger'' than the one constructed in \cite{Kernels_nodal} and \cite{categorical_absorptions}, since we split off one fewer exceptional object than they do.
\end{remark}
\begin{proof}
	Let $n=\dim(X)-1$, for some integer $n \geq 1$. In Lemma \ref{blow_up} we recalled that we can resolve the singularity of $X$ by a single blow-up at the singular point $x \in X$. Let $\pi\colon\widetilde{X}:=\bl_x(X) \to X$ be the corresponding proper birational morphism. We obtain a cartesian diagram
	\begin{equation*}
		\begin{tikzcd}
			Y \arrow[d, "p"'] \arrow[r, hook, "j"] & \widetilde{X}  \arrow[d, "\pi"] \\
			\{x\} \arrow[r, hook, "i"]                     & X,                             
		\end{tikzcd}
	\end{equation*}
	where $Y$ denotes the exceptional divisor of $\widetilde{X}$ and $i$ and $j$ the inclusions of the point $\{x\} \subset X$ and $Y \subset \widetilde{X}$, respectively. In the case that $x\in X$ is an $A_1$ singularity, $Y$ is a smooth quadric, and in the case $x \in X$ is an $A_2$ singularity, $Y$ is a nodal quadric.  Then Proposition \ref{Kuzcomp_for_gorenstein} implies that we have a semiorthogonal decomposition
	\begin{equation}\label{decomp_exceptional}
		D^b(Y) = \langle \mathcal{A}_Y, \mathcal{O}_Y, \mathcal{O}_Y(1), \dots, \mathcal{O}_Y(n-1) \rangle. 
	\end{equation}
	Since $Y$ is Gorenstein in both cases, this decomposition is admissible and we can permute a single component by tensoring with the canonical sheaf $\omega_Y$, by Lemma \ref{SOD_and_Serre}. The adjunction formula yields $\omega_Y = \mathcal{O}_Y(-n)$ and a successive application of the functor $- \otimes \mathcal{O}_Y(-n)$ to the line bundles $\mathcal{O}_Y(1), \dots, \mathcal{O}_Y(n-1)$ gives rise to a semiorthogonal decomposition
	\begin{equation}\label{dual_lef}
		D^b(Y) = \langle \mathcal{O}_Y(1-n), \mathcal{O}_Y(2-n), \dots, \mathcal{O}_Y(-1), \mathcal{A}_Y, \mathcal{O}_Y \rangle.
	\end{equation}
	This is a dual Lefschetz decomposition with respect to the  conormal bundle $\mathcal{N}_{Y/\widetilde{X}}^{\vee} \cong \mathcal{O}_Y(1)$, by setting $\mathcal{B}_{n-1} = \mathcal{B}_{1} = \langle \mathcal{O}_Y\rangle$ and $\mathcal{B}_0 = \langle \mathcal{A}_Y, \mathcal{O}_Y \rangle$. We apply Theorem \ref{catres_withLefschetz} and obtain a semiorthogonal decomposition 
	\begin{equation}\label{SODLefschetz}
		D^b(\widetilde{X})= \langle j_*\mathcal{O}_Y(1-n),\dots, j_*\mathcal{O}_Y(-1),\widetilde{\mathcal{D}}\rangle,
	\end{equation} 
	where 
	\begin{equation*}
		\widetilde{\mathcal{D}} = \{F \in D^b(\widetilde{X}) | \ j^*F \in \mathcal{B}_0\}.
	\end{equation*}
    Let $\pi_*\colon \widetilde{\mathcal{D}} \to \dbx$ denote the restriction of the pushforward functor along $\pi$ to $\widetilde{\mathcal{D}}$.
    Note that we have an inclusion 
    \begin{equation}\label{inclusion}
		p^*(D^{\lilperf}(x)) \subset \mathcal{B}_0,
	\end{equation}
    since all free $\mathcal{O}_Y$-modules are contained in $\mathcal{B}_0$ and for any finite dimensional $k$-vector space $V$ we have $p^*(V) \cong \mathcal{O}_Y^{\oplus \dim(V)}$. Therefore, we can apply Lemma \ref{small_to_big} to obtain an inclusion
    \begin{equation*}
        \im(\pi^* \colon D^{\lilperf}(X) \to D^b(\widetilde{X})) \subset \widetilde{\mathcal{D}}.
    \end{equation*}
    This implies that the triple $(\widetilde{\mathcal{D}}, \pi_*, \pi^*)$ defines a categorical resolution of $\dbx$ by Theorem \ref{catres_withLefschetz}. This resolution is in fact crepant which follows by Proposition \ref{crepancy_check}. Indeed, we have $p^*(D^{\lilperf}(x)) \subset \mathcal{B}_{n-1} \subset \mathcal{B}_0$ and the canonical bundle of $\widetilde{X}$ is given by $\omega_{\widetilde{X}} = \pi^*\omega_X \otimes \mathcal{O}_{\widetilde{X}}((n-1)Y)$, see Lemma \ref{blow_up}.
\end{proof}
The following theorem is the key to determining generators of $\ker(\pi_*) \subset \widetilde{\mathcal{D}}$.
\begin{theorem}[{\cite[Theorem 8.22]{Efimov_key}, \cite[Theorem 5.2]{categorical_absorptions}}]\label{efimov}
	Let $\pi\colon \widetilde{X} \to X$ be a proper birational morphism and $i\colon Z \to X$ a closed immersion such that the schematic preimage $E := \pi^{-1}(Z)$ is a Cartier divisor. Assume that the restriction $\pi\colon \widetilde{X}\setminus E \to X\setminus Z$ is an isomorphism and 
	\begin{equation}\label{good_loci}
		\pi_*\mathcal{O}_{\widetilde{X}}(-mE) \cong \mathcal{J}_Z^m, \quad \text{for all } m \geq 0. 
	\end{equation}
	Consider the cartesian diagram
	\begin{equation*}
		\begin{tikzcd}
			E \arrow[r, "j"] \arrow[d, "p"'] & \widetilde{X} \arrow[d, "\pi"] \\
			Z \arrow[r, "i"]                         & X.               
		\end{tikzcd}
	\end{equation*}
	If the functor $p_*\colon D^b(E) \to D^b(Z)$ is a Verdier localization, then the functor $\pi_*\colon D^b(\widetilde{X}) \to \dbx$ is also a Verdier localization and $\ker(\pi_*)$ is generated by $j_*(\ker(p_*))$. 
\end{theorem}
We will split Theorem \ref{MAIN} into two statements: first we explicitly give a set of the generators of $\ker(\pi_*)$, then we show that they are spherical in Theorem \ref{theorem_spherical} below.
\begin{theorem}\label{theorem_general}
	Let $X$ be a projective variety with an isolated $A_2$ singularity at a point $x \in X$ and let $\dim(X) \geq 2$. Then the crepant categorical resolution $\pi_*\colon \widetilde{\mathcal{D}} \to \dbx$, constructed in Theorem \ref{cuspcatres}, is a Verdier localization. Furthermore, we have:
	\begin{equation}
		\ker(\pi_*) = \begin{cases}
			\langle j_*\mathcal{S} \rangle, \hspace{1.45cm} \text{if } \dim(X) \ \text{is odd,} \\
			\langle j_*\mathcal{S}_1,j_*\mathcal{S}_2 \rangle, \quad \text{if } \dim(X) \ \text{is even,}
		\end{cases}
	\end{equation}
    where $j$ denotes the closed embedding of the exceptional divisor $Y$ into the blow-up of $X$ at $x$.
\end{theorem}
\begin{proof}
    In Lemma \ref{blow_up} we proved that an $A_2$ singularity can be resolved by a single blow-up. Using the notation $\widetilde{X} := \bl_x(X)$, we have a cartesian diagram
    \begin{equation*}
	\begin{tikzcd}
		Y \arrow[r, "j"] \arrow[d, "p"'] & \widetilde{X} \arrow[d, "\pi"] \\
		\{x\} \arrow[r, "i"]                         & X,       
	\end{tikzcd}
    \end{equation*}
    where the exceptional divisor $Y\subset X$ is a nodal quadric. A straightforward computation shows that it suffices to prove that
    \begin{enumerate}\label{conditions_pns}
		\item the canonical map $\mathfrak{m}_{X,x}^m/\mathfrak{m}_{X,x}^{m+1} \to H^0(Y,\mathcal{O}_Y(m))$ is an isomorphism;
		\item $H^i(Y, \mathcal{O}_Y(m)) = 0$ for all $i>0$,
	\end{enumerate}
    for all $m\geq 0$ to show that the requirement (\ref{good_loci}) is satisfied, see \cite[Lemma 5.5]{categorical_absorptions}. These two conditions can be shown analogously to the $A_1$ case, for which we refer to \cite[Lemma 5.7]{categorical_absorptions}. Let $\dim(Y) = n$ and let us restrict $X$ to a formal neighborhood of the $A_2$ singularity $x$, where $X$ can be defined by the (affine) equation $F = x_1^2+\dots+x_{n+1}^2+x_{n+2}^{3} =: q +x_{n+2}^3$. The exceptional divisor $Y = V(q) \subset \mathbb{P}^{n+1}$ induces a long exact sequence of the form
	\begin{equation*}
		\begin{tikzcd}[column sep = 0.5cm, row sep = 0.5cm]
        0  \arrow[r] & {H^0(\mathbb{P}^{n+1},\mathcal{O}(m-2))} \arrow[r] & {H^0(\mathbb{P}^{n+1}, \mathcal{O}(m))} \arrow[r] & {H^0(Y, \mathcal{O}_Y(m))} \arrow[out=-2, in=178, looseness=1.3, overlay]{dll} \\
             & {H^1(\mathbb{P}^{n+1}, \mathcal{O}(m-2)) = 0}      &                                                   &              
        \end{tikzcd}
	\end{equation*}
	for all $m \geq 0$ and the last term of the sequence vanishes for all $m\geq 0$, since $n = \dim(X) -1 \geq 1$. Moreover, there is an isomorphism $H^0(\mathbb{P}^{n+1}, \mathcal{O}(m)) \cong k[x_1,\dots,x_{n+2}]_m$ which, together with the exact sequence, yields 
	\begin{equation*}
		H^0(Y, \mathcal{O}_Y(m)) \cong k[x_1,\dots,x_{n+2}]_m\big{/}q \cdot k[x_1,\dots,x_{n+2}]_{m-2}.
	\end{equation*}
	A simple calculation shows that $\mathfrak{m}_{X,x}^m/\mathfrak{m}_{X,x}^{m+1}$ is canonically isomorphic to the latter quotient.
	The second condition is satisfied since $Y \subset \mathbb{P}^{n+1}$ is a quadric hypersurface. 

    In the following we will assume that $X$ is even dimensional (and the odd dimensional case can be shown analogously). We note that the structure sheaf $\mathcal{O}_Y$ is an exceptional object in the category $D^b(Y)$, so we have a semiorthogonal decomposition $D^b(Y) = \langle \mathcal{O}_Y^\perp, \mathcal{O}_Y \rangle$. This implies that $p_*\colon D^b(Y) \to D^b(\{x\})$ is in fact a Verdier localization with kernel $\ker(p_*) =\langle \mathcal{O}_Y\rangle^\perp$. By an application of Theorem \ref{efimov} we obtain that $\pi_*\colon D^b(\widetilde{X}) \to \dbx$ is a Verdier localization and its kernel is generated by $j_*( \langle \mathcal{O}_Y\rangle^\bot )$. More precisely, by Proposition \ref{Kuznetsov_ff_SOD} and Proposition \ref{Mod_sheaf_equivalence}, there exists a semiorthogonal decomposition 
    \begin{equation*}
		D^b(Y) = \langle \mathcal{O}_Y(1-n), \mathcal{O}_Y(2-n), \dots, \mathcal{O}_Y(-1),\langle \mathcal{S}_1,\mathcal{S}_2 \rangle, \mathcal{O}_Y \rangle
	\end{equation*}
   and with the discussion above we obtain 
    \begin{equation}\label{kernel}
    \begin{aligned}
        \ker(\pi_*) & = j_*\langle \mathcal{O}_Y(1-n), \dots, \mathcal{O}_Y(-1), \mathcal{S}_1, \mathcal{S}_2 \rangle. \\
         &\overset{(*)}{=}\langle j_*\mathcal{O}_Y(1-n),\dots, j_*\mathcal{O}_Y(-1), \langle j_*\mathcal{S}_1, j_*\mathcal{S}_2 \rangle \rangle, 
    \end{aligned}
    \end{equation}
    where $(*)$ holds, because the functor $j_*$ is fully faithful on each of the subcategories $\langle \mathcal{O}_Y(1-n)\rangle , \dots, \langle\mathcal{O}_Y(-1)\rangle$, respectively, see (\ref{SODLefschetz}). In particular, we have $$\widetilde{\mathcal{D}}^\bot \subseteq \ker(\pi_*)$$
    which implies that the restriction $\pi_*\colon\widetilde{\mathcal{D}} \to \dbx$ to the crepant resolution $\widetilde{\mathcal{D}}$ of $\dbx$ is also a Verdier localization. Moreover, there exists a semiorthogonal decomposition 
    \begin{equation*}
        \ker(\pi_*) = \langle \widetilde{\mathcal{D}}^\perp, \ker(\pi_*) \cap \widetilde{\mathcal{D}} \rangle 
    \end{equation*}
    and the comparison with (\ref{kernel}) provides the equality $\ker(\pi_*) \cap \widetilde{\mathcal{D}} = \langle j_*\mathcal{S}_1, j_*\mathcal{S}_2 \rangle$ which is exactly the kernel of the resolution $\pi_*\colon \widetilde{\mathcal{D}} \to \dbx$.
\end{proof}
Before restating the second half of Theorem \ref{MAIN}, we recall the definition of an \textit{$l$-spherical object}. Let $\mathcal{D}$ be a full admissible subcategory of $\dbx$ for a smooth projective variety $X$. Since $\dbx$ admits a Serre functor, it is easy to see that $\mathcal{D}$ also admits a Serre functor $\mathbb{S}_{\mathcal{D}}$, see \cite[Lemma 7.1.14]{geocubic}.
\begin{definition}\label{defi_spherical}
	Let $l \in \mathbb{Z}$. We say that an object $E \in \mathcal{D}$ is $l$-spherical if it satisfies the following properties. 
	\begin{enumerate}
		\item We have an isomorphism 
		\begin{equation*}
			\Hom(E,E[i]) \cong \begin{cases}
				k, \quad \text{if } i = 0,l \\
				0, \quad \hspace{0.03cm}\text{else}.
			\end{cases}
		\end{equation*}
		\item The Serre functor applied to $E$ is given by a shift by $l$, in other words we have 
			$\mathbb{S}_{\mathcal{D}}(E) \cong E[l]$. 
	\end{enumerate}
\end{definition}
\begin{theorem}\label{theorem_spherical}
	In the setting of Theorem \ref{theorem_general}, assume that $X$ is even dimensional. Then the sheaves $j_*\mathcal{S}_1$ and $j_*\mathcal{S}_2$ which generate the kernel $\ker(\pi_*)$, are $2$-spherical.
\end{theorem}
\begin{proof}
    We first show that there exists an isomorphism $\Hom(j_*\mathcal{S}_1, j_*\mathcal{S}_1) \cong k \oplus k[-2]$. Consider the exact triangle
	\begin{equation}\label{triangle_linkage}
		j^*j_*\mathcal{S}_1 \longrightarrow \mathcal{S}_1 \overset{\varepsilon}{\longrightarrow} \mathcal{S}_1\otimes \mathcal{N}_{Y/\widetilde{X}}^\vee[2] = \mathcal{S}_1(1)[2]
	\end{equation}
	which exists since $j \colon Y \to \widetilde{X}$ is a divisorial embedding, see \cite[§3]{exact_triangle}.
    An application of the functor $\Hom(-, \mathcal{S}_1)$ yields 
	\begin{equation}\label{ext_triangle}
		\mathrm{Ext}^{\bullet}(\mathcal{S}_1(1), \mathcal{S}_1)[-2] \overset{\varepsilon^*}{\longrightarrow} \mathrm{Ext}^{\bullet}(\mathcal{S}_1, \mathcal{S}_1) \longrightarrow \mathrm{Ext}^{\bullet}(j^*j_*\mathcal{S}_1, \mathcal{S}_1)
	\end{equation} 
    and by Theorem \ref{heart_intro} we obtain an isomorphism of k-algebras 
    \begin{equation*}
       \mathrm{Ext}^{\bullet}(\mathcal{S}_1, \mathcal{S}_1) \cong k[\theta],
    \end{equation*}
    where the element $\theta$ has degree $2$. 
    The complex $\mathrm{Ext}^{\bullet}(\mathcal{S}_1(1), \mathcal{S}_1)$ admits a structure of a free $k[\theta]$-module of rank $1$, generated by an element in $\mathrm{Ext}^2(\mathcal{S}_1(1), \mathcal{S}_1)$, as it was shown in Lemma \ref{module_structure_twist}. We consider the morphism of $k[\theta]$-modules
	\begin{equation}\label{composition_ext}
		 \varepsilon^* \colon \mathrm{Ext}^{\bullet}(\mathcal{S}_1(1)[2], \mathcal{S}_1) \longrightarrow \mathrm{Ext}^{\bullet}(\mathcal{S}_1, \mathcal{S}_1), \quad f \mapsto f \circ \varepsilon
	\end{equation}
    which is either zero or an isomorphism in degrees $\geq 4$. Suppose that $\varepsilon^*$ is zero. Using the adjunction $j^* \vdash j_*$ we have an isomorphism 
	\begin{equation*}
		\mathrm{Ext}^{\bullet}(j^*j_*\mathcal{S}_1, \mathcal{S}_1) \cong \mathrm{Ext}^{\bullet}(j_*\mathcal{S}_1, j_*\mathcal{S}_1)
	\end{equation*}
	which implies that the complex $\mathrm{Ext}^{\bullet}(j_*\mathcal{S}_1, j_*\mathcal{S}_1)$ in $D^b(\widetilde{X})$ is unbounded. Since $\widetilde{X}$ is smooth and proper, this is a contradiction. This argument works analogously for the sheaf $\mathcal{S}_2$. 
    
    We now show that the second part of Definition \ref{defi_spherical} holds for the sheaf $j_*\mathcal{S}_1$ (resp.\ $j_*\mathcal{S}_2$) and $l=2$.
	We proceed in a similar way as \cite[Lemma 5.10(iii)]{categorical_absorptions} in the case of $A_1$ singularities. 
	For the following computation, let $n = \dim(Y) = \dim(X)-1$ and recall the semiorthogonal decomposition 
    \begin{equation}
		D^b(\widetilde{X})= \langle j_*\mathcal{O}_Y(1-n),\dots, j_*\mathcal{O}_Y(-1),\widetilde{\mathcal{D}}\rangle
	\end{equation}
    which was established in the proof of Theorem \ref{cuspcatres}. We compute 
	\begin{equation}\label{Serre_comp}
		\begin{aligned}
			\mathbb{S}_{\widetilde{\mathcal{D}}}(j_*\mathcal{S}_1) = \mathbb{R}_{\widetilde{\mathcal{D}}^\perp}(\mathbb{S}_{\widetilde{X}}(j_*\mathcal{S}_1)) & = \mathbb{R}_{\widetilde{\mathcal{D}}^\perp} (j_*\mathcal{S}_1 \otimes \omega_{\widetilde{X}}[n+1]) \\
			& = \mathbb{R}_{\widetilde{\mathcal{D}}^\perp}(j_*(\mathcal{S}_1 \otimes j^*\omega_{\widetilde{X}}))[n+1] \\
			& = \mathbb{R}_{\widetilde{\mathcal{D}}^\perp}(j_*\mathcal{S}_1(1-n))[n+1],
		\end{aligned}
	\end{equation}
	where we applied the adjunction formula  
	\begin{equation*}
		j^*\omega_{\widetilde{X}} = \omega_{Y} \otimes j^*\mathcal{O}_{\widetilde{X}}(-Y) = \mathcal{O}_Y(-n) \otimes \mathcal{O}_Y(1) = \mathcal{O}(1-n)
	\end{equation*}
	in the last step. By twisting the short exact sequences (\ref{ses_with_twist}) by $\mathcal{O}_Y(k)$ for a suitable $k\in \mathbb{Z}$ and pushing them forward along $j_*$, we obtain a sequence of morphisms
	\begin{equation}\label{composition}
		j_*\mathcal{S}_1[2] \longrightarrow j_*\mathcal{S}_2(-1)[3] \longrightarrow \dots \longrightarrow j_*\mathcal{S}_1(1-n)[n+1]
	\end{equation} 
	with cones $j_*\mathcal{O}_Y^M(-1)[3], \dots,j_*\mathcal{O}_Y^M(1-n)[n+1]$, respectively. For the cone of the composition of the maps in the above sequence, we have: 
	\begin{equation*}
		\cone(j_*\mathcal{S}_1[2] \longrightarrow j_*\mathcal{S}_1(1-n)[n+1]) \in \langle j_*\mathcal{O}_Y(1-n),\dots, j_*\mathcal{O}_Y(-1) \rangle = \widetilde{\mathcal{D}}^\perp.
	\end{equation*}
	In the proof of Theorem \ref{theorem_general}, we showed that the sheaf  $j_*\mathcal{S}_1$ is an object of the resolution $\widetilde{\mathcal{D}}$. Therefore, we can apply the right mutation functor $\mathbb{R}_{\widetilde{\mathcal{D}}^\perp}$ to the composition of the sequence of maps (\ref{composition}) to obtain an isomorphism 
	\begin{equation*}
		j_*\mathcal{S}_1[2] \overset{\sim}{\longrightarrow} \mathbb{R}_{\widetilde{\mathcal{D}}^\perp}(j_*\mathcal{S}_1(1-n))[n+1]
	\end{equation*}
	which yields  $\mathbb{S}_{\widetilde{\mathcal{D}}}(j_*\mathcal{S}_1) =j_*\mathcal{S}_1[2]$. Since the sequences (\ref{ses_with_twist}) are symmetric in $\mathcal{S}_1$ and $\mathcal{S}_2$, the same argument works for the sheaf $j_*\mathcal{S}_2$. 
	\end{proof}
 \begin{proposition}
     In the setting of Theorem \ref{theorem_general}, assume that $X$ is odd dimensional. Then we have 
     \begin{equation*}
         \mathrm{Ext}^{\bullet}(j_*\mathcal{S}, j_*\mathcal{S}) \cong k \oplus k[-1] \oplus k[-2]
     \end{equation*}
     and $\mathbb{S}_{\widetilde{\mathcal{D}}}(j_*\mathcal{S}) = j_*\mathcal{S}[2]$. In particular, the sheaf $j_*\mathcal{S}$ that generates the kernel $\ker(\pi_*)$ is not $l$-spherical for any natural number $l \geq 0$.
 \end{proposition}
\begin{proof}
    We proceed analogously to the proof of Theorem \ref{theorem_spherical}. There exists an exact triangle 
    \begin{equation*}
		\mathrm{Ext}^{\bullet}(\mathcal{S}(1), \mathcal{S})[-2] \overset{\varepsilon^*}{\longrightarrow} \mathrm{Ext}^{\bullet}(\mathcal{S}, \mathcal{S}) \longrightarrow \mathrm{Ext}^{\bullet}(j^*j_*\mathcal{S}, \mathcal{S}) 
	\end{equation*} 
    and by Theorem \ref{heart_intro} we have an isomorphism of $k$-algebras
    \begin{equation*}
       \mathrm{Ext}^{\bullet}(\mathcal{S}, \mathcal{S}) \cong k[\theta'],
    \end{equation*}
    where the element $\theta'$ has degree $1$. Moreover, by Lemma \ref{even_module_structure_twist} the complex $\mathrm{Ext}^{\bullet}(\mathcal{S}(1), \mathcal{S})$ admits a structure of a free $k[\theta']$-module of rank $1$, generated by an element in $\mathrm{Ext}^1(\mathcal{S}(1), \mathcal{S})$. Therefore it suffices to show that the morphism of $k[\theta']$-modules 
    \begin{equation}\label{even_composition_ext}
		 \varepsilon^* \colon \mathrm{Ext}^{\bullet}(\mathcal{S}(1)[2], \mathcal{S}) \longrightarrow \mathrm{Ext}^{\bullet}(\mathcal{S}, \mathcal{S}), \quad f \mapsto f \circ \varepsilon
	\end{equation}
    is an isomorphism in degrees $\geq 3$. This follows from the same arguments as in the proof of Theorem  \ref{theorem_general}.

    The computation of $\mathbb{S}_{\widetilde{\mathcal{D}}}(j_*\mathcal{S})$ works analogously to the even dimensional case in Theorem \ref{theorem_spherical}.
\end{proof}

\section{Special case of cubic fourfolds}\label{kernels_cubic} 
In this section we apply the main theorem to the case of cubic fourfolds with the aim of proving Theorem \ref{cubic_fourfold_main_result}. The proof consists of two steps: first, we show that the categorical resolution of Theorem \ref{MAIN} restricts to a categorical resolution of the Kuznetsov component of the cubic fourfold and that the resolution is equivalent to the bounded derived category of a K3 surface $S$. Secondly, we will describe the generators of the kernel inside of $D^b(S)$. This generalizes the analogous results in the case of $A_1$ singularities proved in \cite[§5]{dervied_cat_cubic_fourfold} and \cite[§4]{Kernels_nodal}, respectively.

Let $V$ be a $6$-dimensional $k$-vector space and let $X \subset \mathbb{P}(V)$ denote a cubic fourfold with an isolated $A_2$ singularity at a point $x\in X$. We choose projective coordinates $x_0,\dots, x_5$, such that $x = [1:0:\dots:0]$. Let $\pi\colon\widetilde{X} := \bl_x(X) \to X$ denote the corresponding map of the blow-up of $X$ at $x$ which is a resolution of singularities for $X$. Let $\sigma \colon \widetilde{X} \to \mathbb{P}^4$ be the extension of the projection away from the cuspidal point $x$ to the blow-up $\widetilde{X}$. Recall from Proposition \ref{fundamental} that $X$ has a defining equation of the form $x_0q+g$ for a nodal quadric $V(q)$ and some cubic $V(g)$ in $V(x_0) \cong \mathbb{P}^4$, and the complete intersection $S = V(q, g) \subset \mathbb{P}^4$ is a smooth K3 surface by Corollary \ref{S_is_K3}.  

The following lemma was proved in \cite[Lemma 5.1]{dervied_cat_cubic_fourfold} for a cubic fourfold $X$ with a single isolated $A_1$ singularity. 
\begin{lemma}\label{associatedK3}
	The morphism $\sigma$ is isomorphic to the blow-up of $\mathbb{P}^4$ along the K3 surface $S$. We denote by $Y$ and $D$ the exceptional divisors of $\pi$ and $\sigma$ and the corresponding closed immersions by $j\colon Y \hookrightarrow \widetilde{X}$ and $\eta\colon D \hookrightarrow \widetilde{X}$, respectively. There exist two cartesian diagrams
	\begin{equation*}
		\begin{tikzcd}
			& Y \arrow[ld, "p"'] \arrow[r, hook, "j"] & \widetilde{X} \arrow[ld, "\pi"'] \arrow[rd, "\sigma"] & D \arrow[l, hook', "\eta"'] \arrow[rd, "q"] &             \\
			x \arrow[r, hook, "i"] & X                                &                                                       & \mathbb{P}^4                 & S. \arrow[l,hook']
		\end{tikzcd}
	\end{equation*}
The morphism $\sigma \circ j$ identifies $Y$ with the quadric passing through $S$. Moreover, let $H$ and $h$ be pullbacks of classes of hyperplanes in $\mathbb{P}(V)$ and $\mathbb{P}^4$, respectively. Then we have the following relations in $\text{Pic}(\widetilde{X})$:
\begin{equation*}
	Y = 2h - D, \ H=3h-D,\ h=H-Y, \ D=2H-3Y, \ K_{\widetilde{X}} = -5h +D = -3H +2Y.
\end{equation*}
\end{lemma}
\begin{proof}
	For the first two claims, we refer to \cite[Lemma 5.1]{dervied_cat_cubic_fourfold} for a detailed proof. This reference provides a proof in the case of a nodal variety $X$ which works analogously in the cuspidal case, essentially because both singularities can be resolved by one blow-up and they both have multiplicity $2$. For the relations in the Picard group $\picard(\widetilde{X})$, we first note that the right blow-up diagram yields $K_{\widetilde{X}} = -5h + D$. The other relation containing the canonical divisor $K_{\widetilde{X}}$ follows from Lemma \ref{blow_up}, where we proved that the discrepancy of the cusp is $\dim(X)-2=2$.
    Since $\sigma$ is the extension of the projection away from the cuspidal point $x \in X$ to the blow-up $\widetilde{X}$, the relation $h = H - Y$ holds.  
	Finally, note that the proper transform (with respect to $\sigma$) of the quadric $\sigma \circ j(Y) \subset \mathbb{P}^4$ passing through $S$ is contracted by $\pi$. Therefore we have an equation $Y = 2h-D$. The other relations follow by substitution from those we have proved.
\end{proof}
By Proposition \ref{Kuzcomp_for_gorenstein} and Lemma \ref{SOD_and_Serre} there exists an admissible semiorthogonal decomposition
\begin{equation}\label{SODstandard}
	\dbx = \langle \mathcal{A}_X, \mathcal{O}_X,\mathcal{O}_X(H),\mathcal{O}_X(2H) \rangle,
\end{equation}
where the subcategory $\mathcal{A}_X$ denotes the Kuznetsov component of $X$. 
\begin{proof}[Proof of Theorem \ref{cubic_fourfold_main_result}]
    Since the relations in $\picard(\widetilde{X})$ in Lemma \ref{associatedK3} are central for this proof and they coincide with the analogous relations in the $A_1$ case, the proof of Theorem \ref{cubic_fourfold_main_result} for $A_1$ singularities, see \cite[Theorem 5.2]{dervied_cat_cubic_fourfold}, generalizes to the $A_2$ case without substantial changes. Below we will provide a sketch of the proof and refer to \cite[Theorem 5.2]{dervied_cat_cubic_fourfold} for details.
    
	By Theorem \ref{cuspcatres}, there exists a crepant categorical resolution $(\widetilde{\mathcal{D}}, \pi_*,\pi^*)$ together with 
    a semiorthogonal decomposition
	\begin{equation}\label{SODLefschetz2}
		D^b(\widetilde{X})= \langle j_*\mathcal{O}_Y(-2h), j_* \mathcal{O}_Y(-h),\widetilde{\mathcal{D}}\rangle.
	\end{equation} 
	Since the functor $\pi^*\colon D^{\lilperf}(X) \to \widetilde{\mathcal{D}}$ is fully faithful, we can use (\ref{SODstandard}) to produce a semiorthogonal decomposition
	\begin{equation}\label{SOD_Dtilde}
		\widetilde{\mathcal{D}} = \langle \widetilde{\mathcal{A}}_X, \mathcal{O}_{\widetilde{X}},\mathcal{O}_{\widetilde{X}}(H),\mathcal{O}_{\widetilde{X}}(2H) \rangle,
	\end{equation}
	where we define 
    \begin{equation*}
        \widetilde{\mathcal{A}}_X :=\langle \mathcal{O}_{\widetilde{X}}, \mathcal{O}_{\widetilde{X}}(H), \mathcal{O}_{\widetilde{X}}(2H) \rangle^\bot.
    \end{equation*}
    This refines (\ref{SODLefschetz2}) in the following way:
	\begin{equation}\label{SOD1}
		D^b(\widetilde{X})= \langle j_*\mathcal{O}_Y(-2h), j_*\mathcal{O}_Y(-h),\widetilde{\mathcal{A}}_X, \mathcal{O}_{\widetilde{X}},\mathcal{O}_{\widetilde{X}}(H),\mathcal{O}_{\widetilde{X}}(2H) \rangle.
	\end{equation}
    By the definition of $\widetilde{\mathcal{A}}_X$ and the fact that $\pi^*\colon \perf(X) \to \widetilde{\mathcal{D}}$ is fully faithful, we have that $\pi_*(\widetilde{\mathcal{A}}_X) \subset \mathcal{A}_X$ and $\pi^*(\mathcal{A}_X^{\lilperf}) \subset \widetilde{\mathcal{A}}_X$, which implies that $(\widetilde{\mathcal{D}}, \pi_*,\pi^*)$ restricts to a crepant categorical resolution $(\widetilde{\mathcal{A}}_X, \pi_*,\pi^*)$ of the Kuznetsov component $\mathcal{A}_X$.
    
    The equivalence $\widetilde{\mathcal{A}}_X \cong D^b(S)$ can be constructed as follows. By an application of Orlov's blow-up formula to $\sigma\colon \widetilde{X} \to \mathbb{P}^4$, we obtain a semiorthogonal decomposition 
  	\begin{equation}\label{SOD2}
  		D^b(\widetilde{X}) = \langle \Phi(D^b(S)), \mathcal{O}_{\widetilde{X}}(-3h), \mathcal{O}_{\widetilde{X}}(-2h),\mathcal{O}_{\widetilde{X}}(-h),\mathcal{O}_{\widetilde{X}},\mathcal{O}_{\widetilde{X}}(h) \rangle,
  	\end{equation}
  	where $\Phi(\mathcal{F}) = \eta_*q^*\mathcal{F}\otimes \mathcal{O}_{\widetilde{X}}(D)$ for $\mathcal{F} \in D^b(S)$. Then one applies a series of mutations, see \cite[Theorem 5.2]{dervied_cat_cubic_fourfold}, to the decomposition (\ref{SOD2}) to obtain a semiorthogonal decomposition 
  	\begin{equation}\label{sod_after_mut}
  		D^b(\widetilde{X})= \langle j_*\mathcal{O}_Y(-2h), j_*\mathcal{O}_Y(-h), \Phi''(D^b(S)), \mathcal{O}_{\widetilde{X}},\mathcal{O}_{\widetilde{X}}(H),\mathcal{O}_{\widetilde{X}}(2H) \rangle,
  	\end{equation}
  	where 
    \begin{equation*}
        \Phi'' = \mathbb{R}_{\mathcal{O}_{\widetilde{X}}(-h)} \circ \mathbb{R}_{\mathcal{O}_{\widetilde{X}}(-2h)} \circ \mathbb{T}_{\mathcal{O}_{\widetilde{X}}(D-2h)} \circ \eta_* \circ q^*.
    \end{equation*} 
    Here, $\mathbb{T}_{\mathcal{O}_{\widetilde{X}}(D-2h)}$ denotes the functor defined by tensoring with the line bundle $\mathcal{O}_{{\widetilde{X}}(D-2h)}$. Comparing the semiorthogonal decomposition (\ref{sod_after_mut}) with (\ref{SOD1}), we conclude that the functor $\Phi''\colon D^b(S) \to D^b(\widetilde{X})$ induces an equivalence of triangulated categories $D^b(S) \cong \widetilde{\mathcal{A}}_X$.

    The second part of the statement can be proven analogously to the case of a cubic fourfold with an isolated $A_1$ singularity, see \cite[§4]{Kernels_nodal}. We will sketch their proof in the following. The authors first explicitly determine the left adjoint $\Psi\colon \widetilde{\mathcal{A}}_X \to D^b(S)$ of the above equivalence $\Phi''\colon D^b(S) \to \widetilde{\mathcal{A}}_X$, which is 
    given by
     \begin{equation*}
        \Psi= q_* \circ \eta^* \circ \mathbb{T}_{\mathcal{O}_{\widetilde{X}}(-3h+D)[1]} \circ \mathbb{L}_{\mathcal{O}_{\widetilde{X}}(3h-D)} \circ \mathbb{L}_{\mathcal{O}_{\widetilde{X}}(4h-D)}.
    \end{equation*}
    This can be verified in the same way as in the $A_1$ case, see \cite[Proposition 4.3]{Kernels_nodal}. Applying the functor $\Psi$ to $j_*\mathcal{S}_i$, where $i \in \{1,2\}$, shows that these sheaves are invariant under the first two mutations and the twist, as proved in \cite[Proposition 4.4]{Kernels_nodal} for $A_1$ singularities. Therefore we have 
    \begin{equation*}
        \Psi(j_*\mathcal{S}_i) = q_* \circ \eta^*(j_*\mathcal{S}_i), \quad i \in \{1,2\}.
    \end{equation*}
    The claim now follows from the commutative diagram 
    \begin{equation*}
        \begin{tikzcd}
        Y \arrow[d, "j", hook]            & S \arrow[l, "t"', hook] \arrow[d, "\iota", hook] \arrow[dd, "\id", bend left=49] \\
        \widetilde{X} \arrow[d, "\sigma"] & D \arrow[d, "q"] \arrow[l, "\eta", hook]                                         \\
        \mathbb{P}^4                      & S, \arrow[l, hook]                                \end{tikzcd}
    \end{equation*}
    where the upper square is cartesian, see \cite[Lemma 4.1]{Kernels_nodal}. In this situation, we can apply \cite[Proposition 2.29]{Kernels_nodal}, which implies that $\iota_*t^* = \eta^*j_*$. Finally, we have 
    \begin{equation*}
        \Psi(j_*\mathcal{S}_i) = q_* \circ \eta^*(j_*\mathcal{S}_i) = q_* \iota_*t^*\mathcal{S}_i = t^*\mathcal{S}_i, \quad i \in \{1,2\}.
    \end{equation*}
\end{proof}

\clearpage
\printbibliography
\end{document}